\documentclass[reqno,12pt]{article}
\usepackage{amsmath}
\usepackage{amssymb}
\usepackage{amsthm}
\usepackage{graphics, xcolor}
\usepackage{geometry}
\begin{document}

\pagestyle{myheadings}

\setcounter{page}{1}

\def\E{{\,{\bf E}\,}}
\def\P{\,{\bf P}\,}
\def\pr{\,{\bf P}\,}
\def\a{\alpha}
\def\be{\beta}
\def\g{\gamma}
\def\GA{\Gamma}
\def\la{\lambda}
\def\LA{\Lambda}
\def\de{\delta}
\def\om{\omega}
\def\DE{\Delta}
\def\cov{\hbox{\rm cov}\,}
\def\p{\varphi}
\def\<{\bbgl <}\def\8{\left<}\def\9{\right>}
\def\>{\bbgr>}
\def\ovln#1{\,{\overline{\!#1}}}
\def\ov{\overline}
\def\L{{\cal L}}
\def\N{\cal N}
\def\e{\varepsilon}
\def\vk{\varkappa}
\def\t{\tau}
\def\th{\theta}
\def\PM{^{\2\prime}}
\def\BT#1{{\cal B}_d(#1)}
\def\BBBT#1{{\cal B}_d^*(#1)}
\def\BD#1#2{{\cal B}_{#2}(#1)}
\def\BB{{\cal B}_d(\tau)}
\def\BBB{{\cal B}_d^*(\tau)}
\def\AT#1{{\cal A}_d(#1)}
\def\AD#1#2{{\cal A}_{#2}(#1)}
\def\A{{\cal A}_d(\tau)}
\def\ATDR#1#2#3{{\cal A}_{d}''(#1,#2,#3)}
\def\ATQZ#1{{\cal A}_d^*(#1,\de,\rho)}
\def\ADQZ#1#2{{\cal A}_{#2}^*(#1,\de,\rho)}
\def\AQZ{{\cal A}_d^*(\tau,\de,\rho)}
\def\ATDRZ#1#2#3{{\cal A}_{d}^*(#1,#2,#3)}
\def\ADRZ#1#2#3#4{{\cal A}_{#4}^*(#1,#2,#3)}
\def\AU#1#2#3#4{{\cal A}_{#4}^{**}(#1,#2,#3)}
\def\AV#1#2#3{\ovln{{\cal A}}_{#3}(#1,#2)}
\def\Duu#1#2#3{{\<{\bf D}_{#1}{#2},{#3}\>}}
\def\Dpuu#1#2#3{{\<{\bf D}\PM_{#1}{#2},{#3}\>}}
\def\Dhuu#1#2#3{{\<{\bf D}_{#1}(h)\,{#2},\,{#3}\>}}
\def\Dhu#1#2{{\<{\bf D}_{#1}(h)\,{#2},\,{#2}\>}}
\def\Dwhu#1#2#3{{\<{\bf D}_{#1}\bbgl(#3\bbgr)\4{#2},\,{#2}\4\>}}
\def\D #1{{\<{\bf D}\4{#1},{#1}\>}}
\def\Dh #1{{\<{\bf D}(h)\4{#1},{#1}\>}}
\def\Dhuv#1#2{{\<{\bf D}(h)\4{#1},{#2}\>}}
\def\si{\sigma}
\def\Duv#1#2{{\<{\bf D}\4{#1},{#2}\>}}
\def\Du#1#2{{\<{\bf D}_{#1}{#2},{#2}\>}}
\def\Dpu#1#2{{\<{\bf D}\PM_{#1}{#2},{#2}\>}}
\def\norm#1{\left\|#1\right\|}
\def\nnorm#1{\bgl\|#1\bgr\|}
\def\nm#1{\bgl|#1\bgr|}
\def\nrm#1{|\2#1\2|}
\def\nnnorm#1{\left\|\2\smash{#1}\2\right\|}
\def\mnorm#1{\left|\2\smash{#1}\2\right|}
\def\snorm#1{\bbgl\|\2#1\2\bbgr\|}
\def\ssqrt{^{1/2}}
\def\msqrt{^{-1/2}}
\def\me{^{-1}}
\def\htau{\norm h\tau}
\def\ttau{\norm t\tau}
\def\bed#1{\bf 1_{{#1}}}
\def\ind#1{{\bf I}\bgl\{{#1}\bgr\}}
\def\eht{\bbgl(\2\e+\htau\2\bbgr)}
\def\weht#1{\bbgl(\2\e+\nnnorm{#1}\t\2\bbgr)}
 \def\leei#1#2{\log\E e^{i\8\right.\!{#1},{#2}\!\left.\9}}
 \def\eei#1#2{\E e^{i\8\right.\!{#1},{#2}\!\left.\9}}
 \def\lee#1#2{\log\E e^{\8\right.\!{#1},{#2}\!\left.\9}}
 \def\ee#1#2{\E e^{\8\right.\!{#1},{#2}\!\left.\9}}
\def\xxxi#1{{\overline\xi}_{#1}(h)}
\def\xh{{\overline\xi}(h)}
\def\oxh#1#2{{\overline\xi_{#1}}(h_{#2})}
\def\ooxh#1#2{{\overline\xi{}_{#1}\PM}(h'_{#2})}
\def\xxi#1{\xi_#1}
\def\eq{{}}
\def\wth{\widetilde h}
\def\wthk{\widetilde h_k}
\def\wthj#1{\widetilde h_{#1}}
\def\tht{{\th\4\htau}}
\def\twhtk{{\th\4\nnnorm{\wthk}\tau}}
\def\twhtj#1{{\th\4\nnnorm{\wthj#1}\tau}}
\def\wtht{{\th\4\wthk\tau}}
\def\ttt{{\th\4\ttau}}
\def\Rd{{\bf R}^d}
\def\Rk{{\bf R}^k}
\def\Rdm{{\bf R}^{d-1}}
\def\Cd{{\bf C}^d}
\def\Bd{{\cal B}_d}
\def\Gd{\cal G_d}
\def\Fd{{\cal F}_d}
\def\xwth{{\overline\xi}(\wth)}
\def\xwthj#1{{\overline\xi_{{#1}}}(\wthj{#1})}
\def\xhpj#1#2{{\overline\xi{}\PM_{{#1}}}(h'_{#2})}
\def\tdsi{{\t\4d^{3/2}\!\big/\si}}
\def\ftdsi{\ffrac{\t\4d^{3/2}}{\si}}
\def\detD{(2\4\pi)^{d/2}\,(\det {\bf D})\ssqrt}
\def\detDj#1{(2\4\pi)^{d/2}\,(\det {\bf D}_{#1})\ssqrt}
\def\wt{\widetilde}
\def\wh{\widehat}
\def\RE{\hbox{Re}}
\def\IM{\hbox{Im}}
\def\Var{\hbox{\rm Var}\,}
\def\jd{{j=1,\dots,d}}
\def\Y#1#2#3{{#1}_{#2}^{[#3]}}
\def\YY#1#2#3{{#1}_{#2}^{(#3)}}
\def\UU#1#2#3{\wt{#1}_{#2}^{(#3)}}
\def\YYY#1#2#3{{#1}_{#2}^{\{#3\}}}

\def\4{\kern1pt}
\def\gm{\4|\4}
\def\bgm{\,\big|\4}
\def\2{\kern.5pt}
\def\6{\phantom0}
\def\bgr#1{\4\hbox{\large$\bigr#1$}}
\def\bgl#1{\hbox{\large$\bigl#1$}\4}
\def\bbgr#1{\2\hbox{\large$\bigr#1$}}
\def\bbgl#1{\hbox{\large$\bigl#1$}\2}
\def\3{\kern-10pt}
\def\5{\kern-6pt}
\def\7#1{_{(#1)}}
\def\cdt{\4\cdot\4}
\def\R#1#2{\cal R_d\bgl(#1,#2\bgr)}
\def\RB#1#2{\cal R_d\bgl(#1,#2\bgr)}
\def\RBM#1#2{\cal R_{d-1}\bgl(#1,#2\bgr)}
\def\half{{}^1\kern-.5pt\!\big/\!\2{}_2}
\def\sfrac#1#2{{}^{#1}\kern-.5pt\!\big/\!\kern.5pt{}_{#2}}
\def\pfrac#1#2{{#1}\big/{#2}}
\def\ffrac#1#2{\4\frac{\,#1\,}{\,#2\,}\4}
\def\xfrac#1#2{{
\small
\hbox{\4$\frac{\,#1\,}{\,#2\,}$\4}}}
\def\zfrac#1#2{\raise.5pt\hbox{\eightpoint$\dfrac{\,#1\,}{\,#2\,}$}}
\def\step#1{\raise9pt\hbox{\eightpoint${#1}$}}
\def\le{\leq}
\def\ge{\geq}
\def\ph{\phantom}
\def\LL#1{\noindent\llapp{{[#1]{\kern10pt }}}\kern-2pt}
\def\ONE{\hbox{1\kern-2.5pt l}}
\def\bR{\bf R}
\def\BR{\bf R}
\def \CG{\cal G}
\def\cit#1{[\2#1\2]}
\def\={: =}
\def\for{\qquad\hbox{for}\quad}

\newtheorem{theorem}{Theorem}[section]
\newtheorem{lemma}{Lemma}[section]
\newtheorem{remark}{Remark}[section]
\newtheorem{corollary}{Corollary}[section]

\author{F.\ G\"otze$^1$\quad and\quad A.\ Yu.\ Zaitsev$^2$}

\title{Multidimensional Hungarian construction\\ for vectors
with almost Gaussian\\ smooth distributions } \markboth{F.\
G\"otze and A.\ Yu.\ Zaitsev}{Hungarian construction for
 almost Gaussian vectors}

\maketitle

 \noindent {\bf Abstract:} A multidimensional version of
the results of Koml\'os, Major and Tusn\'ady for sums of
independent random vectors with finite exponential moments is
obtained in the particular case where the summands have smooth
distributions which are close to Gaussian ones. The bounds
obtained reflect this closeness. Furthermore, the results provide
sufficient conditions for the existence of i.i.d.\  vectors
\,$X_1, X_2,\dots$ \,with given distributions and corresponding
i.i.d.\  Gaussian vectors \,$Y_1, Y_2,\dots$ \,such that, for
given small~\,$\e$,
$$
{\P\Big\{\4{\limsup\limits_{n\to\infty}
\ffrac1{\log n}\Bigl|\,\sum\limits_{j=1}^n X_j-
\sum\limits_{j=1}^n Y_j\,\Bigr|}\le
\e\Big\}=1}.
$$

\medskip


\vspace*{10pt}

\hrule

\section{Introduction}\label{s1}

\footnotetext
{Research supported by the SFB 343 in Bielefeld
by grant INTAS-RFBR 95-0099 and by grant RFBR-DFG 96-01-00096-ge.
\newline\indent
\,\,$^2$Research supported by the SFB 343 in Bielefeld,
 by the Russian Foundation of Basic Research, grant 96-01-00672,
by grant INTAS-RFBR 95-0099 and by  grant RFBR-DFG 96-01-00096-ge.}

The paper is devoted to an improvement
of a multidimensional version of strong approximation results
of Koml\'os, Major and Tusn\'ady
(KMT)
 for sums of independent
random vectors with finite exponential moments
and with smooth distributions which are close to Gaussian
ones.

Let ${\cal F}_d$ be the set of all $d$-dimen\-sional probability
distributions defined on the $\si$-algebra~\4$\Bd$ \,of Borel
subsets of~\,$\Rd$. \,By \,$\wh F(t)$, \,$t\in\Rd$, we denote the
characteristic function of a~distribution~\,${F\in\Fd}$. The
product of measures is understood as their convolution, that is,
${F\4G=F*G}$. \,The distribution and the corresponding covariance
operator of a random vector~$\xi$ will be denoted by~$\L(\xi)$
and~$\cov \xi$ (or~$\cov F$, if~${F=\L(\xi)}$). The symbol~\,${\bf
I}_d$ \,will be used for the identity operator in~\,${\bf R}^d$.
For~\,$b>0$ \,we denote \,$\log^*b=\max\,\bgl\{1,\,\log b\bgr\}$.
\,Writing \,${z\in\Rd}$ \,(resp.\ $\Cd)$, \,we shall use the
representation \,${z=(z_1,\dots,z_d)=z_1\4e_1+\dots+z_d\4e_d}$,
\,where \,${z_j\in{\bf R}^1}$ \,(resp.\ ${\bf C}^1)$ \,and
the~$e_j$ are the standard orthonormal vectors. The scalar product
is denoted by \,${\langle x,y\rangle=x_1\4\ov y_1+\dots+x_d\4\ov
y_d}$. \,We shall use the Euclidean norm \,$\norm z=\langle
z,z\rangle\ssqrt$ \,and the maximum norm \,${|z|=\max\limits_{1\le
j\le k}\,|z_j|}$. \,The symbols \,$c,c_1,c_2,\dots$ \,will be used
for absolute positive constants. The letter~\,$c$ \,may denote
different constants when we do not need to fix their numerical
values. The ends of proofs will be denoted by~\,$\square$.

Let us consider the definition and some useful properties
of classes of distributions \,${\cal A}_d(\tau)\subset\Fd$, \,$\tau\ge0$,
\,introduced in Zaitsev (1986),
see as well Zaitsev (1995, 1996, 1998a).
The class \,${\cal A}_d(\tau)$ \,(with a fixed \,$\tau\ge0$)
consists of
distributions~$F\in\Fd$ \,for which the function
$$
\p(z)=\p(F,z)=\log\int_{{\bf R}^d}e^{\8z,x\9}F\{dx\}\qquad (\p(0)=0)
$$
is defined and analytic for \,$\norm z \tau<1$,
\,$z\in {{\bf C}\2}^d$, \,and
$$
\bgl|d_ud_v^{\22}\,\p(z)\bgr|\4\le\|u\|\4\tau\,\<{\bf D}\,v,v\>
\qquad \hbox{for all}\ \,u,v\in {\bf R}^d \ \,\hbox{and} \
\,\norm z \tau<1,
$$
where \,${\bf D}=\cov F$, \,and the derivative \,$d_u\p$
\,is given by
$$
d_u\p(z)=\lim_{\be\to 0}\,\ffrac{\p(z+\be\4 u)-\p(z)}\be\,.
$$

It is easy to see that
 \,$\t_1<\t_2$ \,implies \,${\AT{\t_1}\subset\AT{\t_2}}$.
\,Moreover,
 the class~\,$\A$
is closed with respect to convolution:
if~\,$F_1,F_2\in\A$, \,then~\,${F_1\4F_2\in\A}$.
\,The class~\,$\AT0$ \,coincides with the class
of all
Gaussian distributions in~\,$\Rd$. \,The following inequality
 can be considered
as an estimate of the stability of this characterization:
if~\,${F\in\A}$, \,$\t>0$, \,then
$$
\pi\bgl(F,\,\Phi(F)\bgr)\le c\4d^2\t\,\log^*(\t\me),
\eqno(1.1)
$$
where \,$\pi(\4\cdot\4,\4\cdot\4)$ \,is the Prokhorov distance
and~\,$\Phi(F)$ \,denotes
 the Gaussian distribution
whose mean and covariance operator are
the same as those of~\,$F$.
Moreover, for all \,$X\in\Bd$
\,and all~\,$\la>0$, \,we have
\begin{eqnarray*}
F\big\{X\big\}&\le&\hbox{\rlap{\hskip7.25cm(1.2)}}
\Phi(F)\big\{X^\la\big\}
+c\4d^2\exp\Big(-\ffrac\la{c\4d^2\4\t}\Big),\\
\Phi(F)\big\{X\big\}&\le&\hbox{\rlap{\hskip7.25cm(1.3)}}
 F\big\{X^\la\big\}
+c\4d^2\exp\Big(-\ffrac\la{c\4d^2\4\t}\Big),
\end{eqnarray*}
where \,${X^\la=\bgl\{y\in\Rd:\inf\limits_{x\in X}\,\nnnorm{x-y}
<\la\bgr\}}$ \,is the $\la$-neighborhood of the set~\,$X$,
see Zaitsev (1986).

The classes \,$\A$ \,are closely connected with
other natural classes of multidimensional
distributions.
In particular, by the definition of \,$\A$,
 \,any distribution~\,$\L(\xi)$ \,from~\,$\A$
\,has finite exponential moments
\,$\E e^{\8h,\xi\9}$, \,for \,$\htau<1$. \,This leads
to exponential estimates for the tails of
distributions (see, e.g., Lemma~\ref {3.3} below). On the other hand,
if \,$\E e^{\8h,\xi\9}<\infty$, \,for \,${h\in A\subset\Rd}$,
\,where \,$A$ \,is a neighborhood of zero, then \,$F=\L(\xi)\in\AT{\t(F)}$
\,with some~\,$\t(F)$ \,depending on~\,$F$ \,only.

\bigskip

Throughout
we assume that \,$\t\ge0$ \,and
\,$\xi_1,\xi_2,\dots$
\,are random vectors
with given distributions \,${\L(\xi_k)\in \A}$
\,such that \,${\E\xi_k=0}$, \,$\cov\xi_k={\bf I}_d$,
 \,${k=1,2,\dots}$. \,The problem is to construct, for a given
\,$n$, \,$1\le n\le\infty$,
\,on a probability space a~sequence of independent
random vectors \,$ X_1,\dots, X_n$
\,and a sequence
of i.i.d.\  Gaussian random vectors
\,$ Y_1,\dots, Y_n$ \,with
\,$\L(X_k)=\L(\xi_k)$, \,$ \E Y_k=0$,
\,$\cov Y_k={\bf I}_d$, \,$k=1,\dots,n$, \,such that,
with large probability,
$$
\DE(n)=\max_{1\le r\le n}
\,\Bigl|\,\sum\limits_{k=1}^r X_k-\sum\limits_{k=1}^r Y_k\,\Bigr|
$$
is as small as possible.

The aim of the paper is to provide sufficient conditions
for the following Assertion~A:
\smallbreak

\noindent{\bf Assertion A.} {\it There exist absolute positive
constants \,$c_1$, $c_2$ \,and \,$c_3$ \,such that, for
\,${\t\4d^{3/2}\le c_1}$, \,there exists a construction with
$$
\E\exp\Bigl(\ffrac{c_2\,\DE(n)}
{d^{3/2}\4\t}\Bigr)
\le \exp\bgl(c_3\4\log^*d
\,\log^*n\bgr).\eqno(1.4)
$$}

Using the exponential Chebyshev inequality, we see that
~(1.4) implies
$$
{}\P\bgl\{\4c_2\,\DE(n)\ge
\t\4d^{3/2}\bgl(c_3\,\log^*d\,\log^*n
+x\bgr)\4\bgr\}
\le e^{-x},\qquad x\ge0.
\eqno(1.5)
$$
Therefore,
Assertion~A can be considered as a generalization of
the classical result of
 KMT (1975,~1976).
 Assertion~A provides a supplement
to an improvement
of a multidimensional KMT-type result of Einmahl~(1989) presented by
Zaitsev~(1995,~1998a)
 which differs from Assertion~A
by the
restriction~\,$\t\ge1$ \,and
by another explicit power-type dependence of the
constants on the dimension~\,$d$.
\,In a particular case, when \,$d=1$ \,and all summands
have a common variance, the result of Zaitsev
 is equivalent to
the main result of~Sakhanenko~(1984), who extended
the KMT construction to the case of non-identically distributed
summands and stated the dependence of constants
on the distributions of
the summands belonging to a subclass of~\,$\AD\t1$.
\,The main difference between Assertion~A and the aforementioned
 results consists in the fact that
in Assertion~A we consider "small"~\,$\t$,
\,$0\le\t\le c_1\4d^{-3/2}$. \,In previous results the constants are
separated from zero by quantities which are larger than some
absolute constants.
In KMT~(1975,~1976) the dependence of the constants on the
distributions is not specified. From the conditions
$( 1)$ and $( 4)$ in
 Sakhanenko~(1984, Section~1), it~follows that
\,$\Var \xi_k\le\la^{-2}$ \,($\la\me$~\,plays in Sakhanenko's
paper the role of~\,$\t$)
and, if~\,$\Var \xi_k=1$, \,then \,$\la\me\ge1$.
\,This corresponds to the restrictions \,$\a\me\ge2$
\,in Einmahl~(1989, conditions~(3.6) and~(4.3)) and \,$\t\ge1$ \,in
Zaitsev (1995,~1998a, Theorem~1).

Note that in Assertion~A we do not require that
the distributions~\,$\L(\xi_k)$
\,are identical but we assume that they have the same covariance
operators, cf.\  Einmahl~(1989)
and Zaitsev (1995,~1998a). A generalization of the results of
Zaitsev (1995,~1998a) and
of the present paper to the case of non-identical covariance
operators appeared recently in the preprint Zaitsev (1998b).

According to~(1.1)--(1.3), the condition~\,${\L(\xi_k)\in \A}$
\,with small~\,$\t$ \,means that \,$\L(\xi_k)$
\,are close to the corresponding Gaussian laws.
It is easy to see that Assertion~A becomes stronger for
small~\,$\t$ \,(see as well Theorem \ref {1.4}  below). Passing to the limit
as \,$\t\to0$, \,we obtain a spectrum of statements with
the trivial limiting case: if \,$\t=0$ \,(and, hence,
~\,$\L(\xi_k)$ \,are Gaussian) we can take~\,$X_k=Y_k$
\,and~\,${\DE(n)=0}$.

We show that
{\it Assertion\/~{\rm  A} is valid under some additional
 smoothness-type restrictions on~\,$\L(\xi_k)$}.
\,The question about the necessity of these conditions remains open.
The case \,$\t\ge1$ \,considered by Zaitsev (1995,~1998a, Theorem~1)
does not need conditions of such kind.
The formulation of our main result---Theorem~\ref {2.1}---includes some additional notation.
In order to show that the conditions of Theorem~\ref {2.1}
can be verified in some concrete simple situations,
 we consider at first
three particular applications---Theorems~\ref{1.1},
\ref {1.2}  and~\ref {1.3}.

\begin{theorem} \label{1.1}
Assume that the distributions~
\,${\L(\xi_k)\in\A}$
\,can be represented in the form
$$
\L(\xi_k)=H_k\4G,\qquad k=1,\dots,n,
$$
where \,$G$ \,is a Gaussian distribution with covariance operator
\,$\cov G=b^2\,{\bf I}_d$ \,with  $b^2$ satisfying \,${b^2\ge
2^{10}\,\t^2\4d^3\,\log^*\ffrac1\t}$. \,Then Assertion\/ {\rm  A}
is valid.
\end{theorem}

The following example deals with a non-convolution family
of distributions approximating a Gaussian distribution
for small~\,$\t$.

\begin{theorem} \label{1.2}
Let \,$\eta$ \,be a random vector
with an absolutely continuous distribution
 and density
$$
p_\t(x)=\ffrac{\big(4+\t^2\nnnorm x^2\big)
\,\exp\big(-\nnnorm x^2\!/2\big)}
{(2\4\pi)^{d/2}\4(4+\t^2\4d)},\qquad x\in\Rd.
\eqno(1.6)
$$
Assume that~
\,${\L(\xi_k)=\L\big(\eta/\gamma\big)}$,
\,$k=1,\dots,n$, \,where
$$
\gamma^2=\ffrac{\big(4+\t^2\4(d+2)\big)}
{(4+\t^2\4d)},\qquad \gamma>0.
\eqno(1.7)
$$
Then Assertion\/ {\rm  A} is valid.
\end{theorem}

The proof of Theorem \ref{1.2} can be apparently extended
to the distributions with some
more generale densities of type \,$P(\t^2\nnnorm x^2)\,
\,\exp\big(-c\,\nnnorm x^2\big)$, \,where \,$P(\cdt)$ \,
is a suitable polynomial.

\begin{theorem} \label{1.3}
Assume that a random vector \,$\zeta$ \,satisfies
the relations
$$
\E\zeta=0,\qquad
\P\bgl\{\nnnorm \zeta\le b_1\bgr\}=1, \qquad H\=\L(\zeta)\in\AT{b_2}
\eqno(1.8)
$$
and admits a differentiable density~\,$p(\cdt)$ \,such that
$$
\sup_{x\in\Rd}\,\bgl|d_u\,p(x)\bgr|\le b_3\,\nnnorm u,\qquad
 \hbox{for all}\quad
u\in\Rd,
\eqno(1.9)
$$
with some positive \,$b_1,\,b_2$ \,and \,$b_3$.
Let \,$ \zeta_1,\zeta_2,\dots$ \,be independent
copies of~\,$\zeta$. Write
$$
\t=b_2\4m\msqrt,
\eqno(1.10)
$$
where
\,$m$ \,is a positive integer.
Assume that the distributions~
\,${\L(\xi_k)}$ \,can be represented in the form
$$
\L(\xi_k)=\YY L{}k\4P,\qquad k=1,\dots,n,\nopagebreak
\eqno(1.11)
$$
where
$$
\YY L{}k\in\A\quad\hbox{and}\quad
P=\L\bgl(\bgl(\zeta_1+\dots+\zeta_m\bgr)\big/\sqrt m\bgr).
\eqno(1.12)
$$
Then there exist a positive
 \,$b_4$ \,depending on~\,$H$
\,only and such that \,$m\ge b_4$
\,implies Assertion\/~{\rm  A}.
\end{theorem}

\begin{remark} \label{r1.1}
If all the distributions~\,$\YY L{}k$
\,are concentrated at zero,
then the statement of Theorem\/~{\rm \ref{1.3}}
$($for~\,$\t=b\4m\msqrt$ \,with some~\,$b=b(H))$
can be derived
from the main results of\/ {\rm  KMT}\/ {\rm (1975,~1976)} $($for \,$d=1)$
and of\/ {\rm  Zaitsev}\/~{\rm (1995,~1998a)} $($for~\,$d\ge1).$
\end{remark}
\smallbreak

A consequence of Assertion~A is given in Theorem~\ref {1.4}
below.

\begin{theorem} \label{1.4}
Assume that \,$\xi,\, \xi_1,\xi_2,\dots,$ \,are i.i.d.\ random
vectors with a common distribution \,${\L(\xi)\in \A}$. \,Let
 Assertion\/~{\rm  A} be satisfied for~
\,$\xi_1,\dots,\xi_{n}$ \,for all~\,$n$ \,with some
\,$c_1$, \,$c_2$ \,and~\,$c_3$ \,independent of~\,$n$.
\,Suppose that \,$\t\4d^{3/2}\le c_1$. \,Then
there exist a construction such that
$$
{}\P\Big\{\4\limsup_{n\to\infty}
\ffrac1{\log n}\Bigl|\,\sum\limits_{j=1}^n X_j-
\sum\limits_{j=1}^n Y_j\,\Bigr|\le
c_4\,\t\4d^{3/2}\log^*d\,\Big\}=1
\eqno(1.13)
$$
with some constant~\,$c_4=c_4(c_2,c_3)$.
\end{theorem}

From a result of B\'artfai~(1966) it follows
that the rate \,$O(\log n)$
\,in~(1.13) is the best possible if
~\,$\L(\xi)$ \,is non-Gaussian. In the case of distributions
with finite exponential moments this rate
was established by Zaitsev~(1995,~1998a, Corollary~1).
Theorems \ref {1.1}--\ref {1.3}  and~\ref {2.1} provide examples of
smooth distributions
which are close to Gaussian ones and for which the constants
corresponding to this rate are arbitrarily small.
The existence of such examples has been already
 mentioned in the one-dimensional case, e.g., by Major~(1978, p.~498).

The paper is organized as follows. In Section~\ref{s2}
we formulate Theorem~\ref {2.1}. To this end
we define at first a class of distributions~\,$\AV \t\rho d$
\,used in Theorem~\ref {2.1}.
The definition of this class is given in terms of
smoothness conditions on the so-called conjugate distributions.
Then we describe a multidimensional version of the KMT dyadic scheme,
cf.\  Einmahl~(1989).
 We prove Theorem~\ref {2.1} in Section~\ref{s3}.
Section~\ref{s4} is devoted to the proofs of Theorems \ref {1.1}--\ref{1.4}.

A preliminary version of the present paper appeared as
the preprint G\"otze and Zaitsev (1997).

\medbreak

\noindent
{\bf Acknowledgment} \
The authors would like
to thank V.~Bentkus for very useful discussions.

\vspace*{30pt}

\section{The main result}\label{s2}

Let \,$F=\L(\xi)\in\A$,
\,$\norm h\t<1$, \,$h\in {\bf R}^d$.
\,The conjugate distribution
\,$\ovln F=\ovln F(h)$ \,is defined by
$$
\ovln F\{dx\}=
\bgl(\!\E e^{\8h,\xi\9}\!\bgr)^{-1}e^{\8h,x\9}F\{dx\}.
\eqno(2.1)
$$
Sometimes
we shall write \,$F_h=\ovln F(h)$.
\,It is clear that \,$\ovln F(0)=F$.
\,Denote by \,${\overline{\xi}}(h)$ \,a random vector with
\,$\L\bbgl(\2\overline{\xi}(h)\bbgr)=\ovln F(h)$.
\,From (2.1) it follows that
$$
\E f\big(\overline{\xi}(h)\big)=\bgl(\!\E e^{\8h,\xi\9}\!\bgr)^{-1}
\E f(\xi)\,e^{\8h,\xi\9},
\eqno(2.2)
$$
provided that \,$\E \bgl| f(\xi)\,e^{\8h,\xi\9}\bgr|<\infty$.
\,It is easy to see that
$$
\hbox{if}\quad U_1,U_2\in\A,\quad U=U_1\4U_2,\quad
\hbox{then}\quad\ovln U(h)=\ovln U_1(h)\,\ovln U_2(h).
\eqno(2.3)
$$

Below we shall also use the following subclasses of~\,$\A$
\,containing distributions satisfying some special smoothness-type
restrictions.
Let \,$\t\ge0$, \,$\de>0$, \,${\rho>0}$, \,$ h\in\Rd$.
\,Consider the conditions:
$$
\int\limits_{\rho\2\norm t\2\t\2d\ge1}\bbgl|\wh F_h(t)\bbgr|
\,dt\le\ffrac{(2\4\pi)^{d/2}\4\t\4d^{3/2}}
{\si\,(\det{\bf D})^{1/2}},\nopagebreak
\eqno(2.4)
$$
$$
\int\limits_{\rho\2\norm t\2\t\2d\ge1}\bbgl|\wh F_h(t)\bbgr|
\,dt\le\ffrac{(2\4\pi)^{d/2}\4\t^2\4d^{2}}{\si^2\,(\det{\bf D})^{1/2}},
\nopagebreak
\eqno(2.5)
$$
$$
\int\limits_{\rho\2\norm t\2\t\2d\ge1} \bbgl|\langle
t,v\rangle\4\wh F_h(t)\bbgr| \,dt\le\ffrac{(2\4\pi)^{d/2}\4\<{\bf
D}\me v,v\>^{\!1/2}} {\de\,(\det{\bf D})^{1/2}}, \qquad\hbox{for
all}\quad v\in\Rd,\nopagebreak \eqno(2.6)
$$
where \,$F_h=\ovln F(h)$ and
 \,$\si^2=\si^2(F)>0$ \,is the minimal eigenvalue
of~\,${\bf D}=\cov F$.
\,Denote by \,$\AV \t\rho d$ \,(resp.\  \,$\AQZ$)
\,the class of distributions \,$F\in\A$
such that the condition~(2.4)
 (resp.\  (2.5) and~(2.6)) is satisfied
 for \,$ h\in\Rd$, \,$\norm h\t<1$. \,It is easy to see that
$$
\AQZ\subset\AV \t\rho d,
\qquad\hbox{if}\quad \ffrac{\t\4d^{1/2}}{\si}\le1.
\eqno(2.7)
$$
In this paper the class \,$\AV \t\rho d$ \,plays the role
of the class \,$\AQZ$ \,which was used by Zaitsev (1995,~1998a),
see also Sakhanenko (1984, inequality~(49),~ p.~9) or
Einmahl (1989, inequality~(1.5)).
Note that~(2.2) implies
$$
\wh F_h(t)=\E e^{\8it,\4\smash{\ov\xi(h)}\9}
=\bgl(\!\E e^{\8h,\4\xi\9}\!\bgr)^{-1}
\E e^{\8h+it,\4\xi\9}.
\eqno(2.8)
$$
\bigskip

\noindent{\bf The dyadic scheme.}
Let \,$N$ be a positive integer and
\,$\bgl\{\xi_1,\dots,\xi_{2^N}\bbgr\}$
\,a collection of $d$-dimensional independent random vectors.
 Denote
$$
\wt S_0=0;\qquad
\wt S_k=\sum_{l=1}^k\xi_l,\quad 1\le k\le2^N;
\eqno(2.9)
$$
$$
U_{m,k}^*=\wt S_{(k+1)\cdot2^m}-\wt S_{k\cdot2^m},
\qquad 0\le k<2^{N-m},\quad0\le m\le N.
\eqno(2.10)
$$
In particular, \,$U_{0,k}^*=\xi_{k+1}$, \,$U_{N,0}^*=\wt
S_{2^N}=\xi_{1}+\dots+\xi_{2^N}$. \,In the sequel we call {\it
block of summands} a collection of summands with indices of the
form \,${k\cdot2^m+1,\dots,(k+1)\cdot2^m}$, \,where \,$0\le
k<2^{N-m}$, \,$0\le m\le N$. \,Thus, \,$U_{m,k}^*$ is the sum over
a block containing $2^m$~summands. Put
$$
\wt U_{n,k}^*
=U_{n-1,2k}^*-U_{n-1,2k+1}^*,
\qquad 0\le k<2^{N-n},\quad1\le n\le N.
\eqno(2.11)
$$
Note that
$$
U_{n-1,2k}^*+U_{n-1,2k+1}^*=U_{n,k}^*,
\qquad 0\le k<2^{N-n},\quad1\le n\le N.
\eqno(2.12)
$$
Introduce the vectors
$$
\wt{\bf U}_{n,k}^*
=\bgl(U_{n-1,2k}^*,\,U_{n-1,2k+1}^*\big)\in{\bf R}^{2d},
\qquad 0\le k<2^{N-n},\quad1\le n\le N,
\eqno(2.13)
$$
with the first \,$d$ \,coordinates
coinciding with those of the vectors
\,$U_{n-1,2k}^*$
\,and with the last \,$d$ \,coordinates
coinciding with those of the vectors
\,$U_{n-1,2k+1}^*$. \,Similarly, denote
$$
{\bf U}_{n,k}^*
=\bgl(U_{n,k}^*,\,\wt U_{n,k}^*\big)\in{\bf R}^{2d},
\qquad 0\le k<2^{N-n},\quad1\le n\le N.
\eqno(2.14)
$$

Introduce now the
projectors \,${\bf P}_{\!i}:{\bf R}^s\to{\bf R}^1$
and \,$\ov{\bf P}_j:{\bf R}^s\to{\bf R}^j$, \,for
\,$i,\,j=1,\dots,s$, \,by the relations
\,${\bf P}_{\!i}\4 x=x_i$,
\,$\ov{\bf P}_j\4 x=(x_ 1,\dots,x_ j)$,
\,\,where
${x=(x_1,\dots,x_s)\in{\bf R}^s}$
\,(we shall use this notation for \,$s=d$ \,or \,$s=2\4d$). \

It is easy to see that, according to~(2.11)--(2.14),
$$
{\bf U}_{n,k}^*={\bf A}\,\wt{\bf U}_{n,k}^*\in{\bf R}^{2d},
\qquad 0\le k<2^{N-n},\quad1\le n\le N,
\eqno(2.15)
$$
where \,${\bf A}:{\bf R}^{2d}\to {\bf R}^{2d}$
\,is a linear operator defined,
for \,$x=(x_1,\dots,x_{2d})\in{\bf R}^{2d}$,
\,as~follows:
$$
\begin{tabular}{r c l}
${\bf P}_{\!j}\,{\bf A}\,x\5$ &=& $\5x_{j}+x_{d+j},\qquad
\jd,$\\
${\bf P}_{\!j}\,{\bf A}\,x\5$ &=& $\5x_{j}-x_{d+j},\qquad
j=d+1,\dots,2\4d.$
\end{tabular}
\eqno(2.16)
$$
Denote
$$
\begin{tabular}{r c l}
${\bf U}_{n,k}^{*(j)}\5$ &=& $\5{\bf P}_j\,{\bf U}_{n,k}^*,$\\
${\bf U}_{n,k}^{*j}\5$ &=& $\5\bgl({\bf U}_{n,k}^{*( 1)},\dots,
{\bf U}_{n,k}^{*(j)}\bgr)=\ov{\bf P}_j\,{\bf U}_{n,k}^*
\in{\bf R}^{j},$
\end{tabular}
\qquad j=1,\dots,2\4d.
\eqno(2.17)
$$

Now we can formulate the main result of the paper.

\begin{theorem} \label{2.1}
Let the conditions described in\/ {\rm (2.9)--(2.17)} be satisfied,
 \,$\tau\ge0$ \,and
\,$\E\xi_k=0$,
\,$\cov \xi_k={\bf I}_d$, \,$k=1,\dots,2^N$.
\,Assume that
$$
{\cal L}\big({\bf U}_{n,k}^{*j}\big)\in \AV\t4j
\for 0\le k<2^{N-n},\quad1\le n\le N,\quad
d\le j\le 2\4d,
\eqno(2.18)
$$
and
$$
{\cal L}\big({\bf U}_{N,0}^{*j}\big)\in \AV\t4j
\for
1\le j\le 2\4d.
\eqno(2.19)
$$
Then there
exist absolute positive constants \,$c_5$, $c_6$ \,and \,$c_7$
\,such that, for \,${\t\4d^{3/2}\le c_5}$, \,one can construct
on a probability space sequences of independent
random vectors \,$ X_1,\dots, X_{2^N}$
\,and
 i.i.d.\  Gaussian random vectors
\,$ Y_1,\dots, Y_{2^N}$ \,so that
$$
\L(X_k)=\L(\xi_k), \quad \E Y_k=0, \quad\cov Y_k={\bf I}_d, \qquad
k=1,\dots,2^N,
\eqno(2.20)
$$
and
$$
\E\exp\Bigl(\ffrac{c_6\,\DE(2^N)}
{d^{3/2}\4\t}\Bigr)
\le \exp\bgl(c_7\4N\,\log^*d\bgr),
\eqno(2.21)
$$
where \,$\DE(2^N)=\max\limits_{1\le r\le 2^N}
\,\Bigl|\,\sum\limits_{k=1}^r X_k-\sum\limits_{k=1}^r Y_k\,\Bigr|$.
\end{theorem}

\bigskip

Theorem~\ref {2.1} says that the conditions~(2.18)
and~(2.19) suffice for Assertion~A.
However, these conditions require
that the number of summands is~\,$2^N$.
\,For an arbitrary
number of summands, one should consider additional
(for simplicity, Gaussian) summands in order to apply Theorem~\ref {2.1}.

Below we shall prove Theorem~\ref {2.1}.
Suppose that its conditions are
satisfied.

At first, we describe a procedure of constructing the random
vectors~$\big\{U_{n,k}\big\}$ \,with distributions
\,$\L\bgl(\big\{U_{n,k}\big\}\bgr)
=\L\bgl(\big\{U_{n,k}^*\big\}\bgr)$, \,provided that the vectors~
\,$Y_1,\dots,Y_{2^N}$ \,are already constructed (then we shall
define \,$X_k=U_{0,k-1}$, \,$k=1,\dots,2^N$). This procedure is an
extension of the KMT~(1975,~1976) dyadic scheme to the
multivariate case due to Einmahl~(1989). For this purpose we shall
use the so-called Rosenblatt quantile transformation
(see~Rosenblatt~(1952) and Einmahl~(1989)).

Denote by \,$F_{N,0}^{( 1)}(x_1)
=\P\bgl\{{\bf P}_1\,U_{N,0}^*<x_1\bgr\}$,
\,$x_1\in{\bf R}^1$,
\,the distribution function of the first
coordinate of the vector~$U_{N,0}^*$.
 Introduce the conditional distributions,
denoting by
\,$ F_{N,0}^{(j)}\bgl(\,\cdot\,\bgm x_1,\dots,x_{j-1}\bgr)$,
\,$2\le j\le d$, \,the regular conditional
distribution function (r.c.d.f.)
of~\,${\bf P}_j\,U_{N,0}^*$, \,given
\,$\ov{\bf P}_{j-1}\,U_{N,0}^*=(x_1,\dots,x_{j-1})$.
\,Denote by
\,$\wt F_{n,k}^{(j)}\bgl(\,\cdot\,\bgm
x_1,\dots,x_{j-1}\bgr)$ \,the r.c.d.f. of~
\,${\bf P}_{\!j}\,{\bf U}_{n,k}^*$,
\,given
\,$\ov{\bf P}_{j-1}\4{\bf U}_{n,k}^*
=(x_1,\dots,x_{j-1})$, \,for~
\,${0\le k<2^{N-n}}$, \,$1\le n\le N$, \,$d+1\le j\le 2\4d$.
\,Put
$$
 T_k=\sum_{l=1}^kY_l,\qquad 1\le k\le2^N;\nopagebreak
\eqno(2.22)
$$
$$
\begin{tabular}{r c l}\vspace{5pt}
$V_{m,k}\5$ &=& $\5\bgl(V_{m,k}^{( 1)},\dots,V_{m,k}^{(d)}\bgr)
=T_{(k+1)\cdot2^m}-T_{k\cdot2^m},$ \\ \vspace{15pt}
& &
\hskip120pt
$0\le k<2^{N-m}\!,\quad0\le m\le N;$\\ \vspace{5pt}
$\wt{\bf V}_{n,k}\5$
&=& $\5\bgl(V_{n-1,2k},\,V_{n-1,2k+1}\big)
=\bgl(\wt{\bf V}_{n,k}^{( 1)},\dots,\wt{\bf V}_{n,k}^{(2d)}\bgr)
\in{\bf R}^{2d},$ \\
& &
\hskip120pt
$0\le k<2^{N-n},\quad1\le n\le N;$
\end{tabular}\nopagebreak
\eqno(2.23)
$$
and
$$
{\bf V}_{n,k}=\bgl({\bf V}_{n,k}^{( 1)},\dots,
{\bf V}_{n,k}^{(2d)}\bgr)={\bf A}\4\wt{\bf V}_{n,k}\in{\bf R}^{2d},
\qquad 0\le k<2^{N-n},\quad1\le n\le N.
\eqno(2.24)
$$
According to the definition of the operator~\,${\bf A}$, \,we have
(see~(2.11)--(2.16) and~(2.22)--(2.24))
$$
{\bf V}_{n,k}
=\bgl(V_{n,k},\,\wt V_{n,k}\big)\in{\bf R}^{2d},
\qquad 0\le k<2^{N-n},\quad1\le n\le N,
\eqno(2.25)
$$
where
$$
\begin{tabular}{r c l}
$V_{n,k}\5$ &=& $\5V_{n-1,2k}+V_{n-1,2k+1},$\\
$\wt V_{n,k}\5$ &=& $\5V_{n-1,2k}-V_{n-1,2k+1},$
\end{tabular}
\qquad 0\le k<2^{N-n},\quad1\le n\le N,
\eqno(2.26)
$$
and
$$
V_{N,0}=Y_1+\dots+Y_{2^N}.
\eqno(2.27)
$$
Thus, the vectors
\,$V_{m,k},\wt{\bf V}_{n,k}$ and ${\bf V}_{n,k}$
\,can be  constructed from the vectors~$\4Y_1,\dots,Y_{2^N}$
\,by the same linear procedure which was used for constructing
the vectors~
\,$U_{m,k}^*,\wt{\bf U}_{n,k}^*$ and ${\bf U}_{n,k}^*$
\,from the vectors \,$\xi_1,\dots,\xi_{2^N}$.

It is obvious that, for fixed \,$n$ and~\,$k$,
$$
\cov {\bf U}_{n,k}^*=\cov {\bf V}_{n,k}=2^n\,{\bf I}_{2d}
\eqno(2.28)
$$
and, hence, the coordinates of the
Gaussian vector~\,${\bf V}_{n,k}$
\,are independent with the same distribution function~
\,$\Phi_{2^{n/2}}(\cdt)$
\,(here and below
$$
\Phi_{\si}(x)=
\int\limits_{-\infty}^{x}\ffrac1{\sqrt{2\4\pi}\,\si}
\,\exp\Bigl(-\ffrac{y^2}{2\4\si^2}\Bigr)\,dy,
\qquad x\in{\bf R}^1,\quad\si>0,
$$
is the distribution function
of the normal law with mean zero and variance~\,$\si^2$).

Denote now the new collection of random vectors~$X_k$
as follows. At first we define
$$
\begin{tabular}{r c l}\vspace{5pt}
$U_{N,0}^{( 1)}\5$ &=& $\5\bgl(\!F_{N,0}^{( 1)}\bgr)^{\!-1}
\bgl(\Phi_{2^{N/2}}\big(V_{N,0}^{( 1)}\big)\bgr)\qquad
\hbox{and, \quad for} \quad 2\le j\le d,$\\
$U_{N,0}^{(j)}\5$ &=& $\5\bgl(\! F_{N,0}^{(j)}\bgr)^{\!-1}
\bgl( \Phi_{2^{N/2}}\big(V_{N,0}^{(j)}\big)\bgr|
\,U_{N,0}^{( 1)},\dots,U_{N,0}^{(j-1)}\bgr)$
\end{tabular}
\eqno(2.29)
$$
(here \,$\bgl(\!F_{N,0}^{( 1)}\bgr)^{\!-1}(t)=
\sup\,\bgl\{x:F_{N,0}^{( 1)}(x)\le t\bgr\}$, \,$0<t<1$,
\,and so on).
Taking into account that the distributions of the random vectors
 \,$\xi_1,\dots,\xi_{2^N}$
\,are absolutely continuous, we see that
this formula can be rewritten in a more natural
 form, cf.\ Sakhanenko (1984,~p.\ 30--31):
$$
\begin{tabular}{r c l}\vspace{5pt}
$F_{N,0}^{( 1)}\big(U_{N,0}^{( 1)}\big)\5$ &=&
$\5\Phi_{2^{N/2}}\big(V_{N,0}^{( 1)}\big),$ \\
 $F_{N,0}^{(j)}\big(U_{N,0}^{(j)}\bgm
U_{N,0}^{( 1)},\dots,U_{N,0}^{(j-1)}\big)\5$ &=&
$\5 \Phi_{2^{N/2}}\big(V_{N,0}^{(j)}\big),\for 2\le j\le d.$
\end{tabular}
\eqno(2.30)
$$

Suppose that the random vectors
$$
U_{n,k}=\bgl(U_{n,k}^{( 1)},\dots,U_{n,k}^{(d)}\bgr),\qquad
0\le k<2^{N-n},
\eqno(2.31)
$$
corresponding to blocks
containing each \,$2^{n}$ \4summands
with fixed~\,$n$, \,${1\le n\le N}$, \,are already constructed.
 Now our aim is to construct
the blocks containing each \,$2^{n-1}$ \4summands.
To this end we define
$$
{\bf U}_{n,k}^{(j)}={\bf P}_{\!j}\,U_{n,k}= U_{n,k}^{(j)},
\quad 1\le j\le d,
\eqno(2.32)
$$
and, for \,$d+1\le j\le 2\4d$,
$$
{\bf U}_{n,k}^{(j)}=\bgl(\!\wt F_{n,k}^{(j)}\bgr)^{\!-1}
\bgl(\Phi_{2^{n/2}}\big({\bf V}_{n,k}^{(j)}
\big)\bgr|\4
{\bf U}_{n,k}^{( 1)},\dots,
{\bf U}_{n,k}^{(j-1)}\bgr) .
\eqno(2.33)
$$
It is clear that~(2.33) can be rewritten in a form similar to~(2.30).
Then we put
$$
\begin{tabular}{r c l}
${\bf U}_{n,k}\5$ &=& $\5\bgl({\bf U}_{n,k}^{( 1)},\dots,
{\bf U}_{n,k}^{(2d)}\bgr)\in{\bf R}^{2d},$ \\
${\bf U}_{n,k}^j\5$ &=& $\5\bgl({\bf U}_{n,k}^{( 1)},\dots,
{\bf U}_{n,k}^{(j)}\bgr)=\ov{\bf P}_j\,{\bf U}_{n,k}
\in{\bf R}^{j},\qquad j=1,\dots,2\4d,$\\
 $\wt U_{n,k}^{(j)}\5$ &=& $\5{\bf U}_{n,k}^{(j+d)},
\qquad j=1,\dots,d,$
\\
$\wt U_{n,k}\5$ &=& $\5\bgl(\wt U_{n,k}^{( 1)},\dots,
\wt U_{n,k}^{(d)}\bgr)\in{\bf R}^{d}$
\end{tabular}
\eqno(2.34)
$$
and
$$
\begin{tabular}{r c l}\vspace{5pt}
$U_{n-1,2k}\5$ &=& $
\5\bgl(U_{n,k}+\wt U_{n,k}\bgr)/2,$\\
$U_{n-1,2k+1}\5$ &=& $\displaystyle
\5\bgl(U_{n,k}-\wt U_{n,k}\bgr)/2.$
\end{tabular}
\eqno(2.35)
$$
Thus, we have constructed the random vectors
\,$U_{n-1,k}$, \,$0\le k<2^{N-n+1}$. \,After \,$N$~\,steps we
 obtain the random vectors~\,$U_{0,k}$, \,$0\le k<2^{N}$.
\,Now we set
$$
X_k=U_{0,k-1},\qquad
S_0=0,\quad S_k=\sum_{l=1}^kX_l,\qquad1\le k\le2^{N}.
\eqno(2.36)
$$

\begin{lemma} \label{l2.1} {\rm (Einmahl~(1989))}
The joint distribution of the vectors \,$U_{n,k}$ \,and~\,${\bf
U}_{n,k}$ \,coincides with that of the vectors~\,$U_{n,k}^*$
\,and~\,${\bf U}_{n,k}^*$. \,In particular, \,$X_k$,
\,$k=1,\dots,2^N$, \,are independent and \,$\L(X_k)=\L(\xi_k)$.
\end{lemma}

Moreover, according to~(2.11) and~(2.12), we have
$$
\begin{tabular}{r c l}
\vspace{5pt}
$\wt U_{n,k}\5$ &=& $\5U_{n-1,2k}-U_{n-1,2k+1},$\\
$U_{n,k}\5$ &=& $\5U_{n-1,2k}+U_{n-1,2k+1}
=S_{(k+1)\cdot2^n}-S_{k\cdot2^n},$
\end{tabular}
\eqno(2.37)
$$
for \,$ 0\le k<2^{N-n}$, \,$1\le n\le N$
\,(it is clear that (2.37) follows from (2.35)).
Furthermore, putting
$$
\wt{\bf U}_{n,k}
=\bgl(U_{n-1,2k},\,U_{n-1,2k+1}\big)\in{\bf R}^{2d},
\qquad 0\le k<2^{N-n},\quad1\le n\le N,\nopagebreak
\eqno(2.38)
$$
we have (see~(2.13) and~(2.15))
$$
{\bf U}_{n,k}={\bf A}\,\wt{\bf U}_{n,k}\in{\bf R}^{2d},
\qquad 0\le k<2^{N-n},\quad1\le n\le N.
\eqno(2.39)
$$
Note that it is not difficult to verify that,
according to~(2.16),
$$
\nnnorm{\bf A}=\ffrac1{\nnnorm{{\bf A}\me}}=
\nnnorm{{\bf A}^*}=\ffrac1{\nnnorm{({\bf A}^*)\me}}=
\sqrt2,
\eqno(2.40)
$$
where the asterisk is used to denote the adjoint operator~\,${\bf A}^*$
for the operator~\,${\bf A}$.

\begin{remark} \label{r2.1}
The conditions of Theorem\/~{\rm \ref {2.1}} imply the coincidence of the
corresponding first and second moments
 of the vectors
\,${\bf U}=\bgl\{U_{n,k},\,\wt{\bf U}_{n,k},\,{\bf U}_{n,k}\bgr\}$
\,and
\,${\bf V}=\bgl\{V_{n,k},\,\wt{\bf V}_{n,k},\,{\bf V}_{n,k}\bgr\}$
\,since the vectors~\,${\bf U}$ \,can be restored from vectors
~\,$X_1,\dots,X_{2^N}$ \,by the same linear procedure
which is used for reconstruction of the vectors~\,${\bf V}$
\,from \,$Y_1,\dots,Y_{2^N}$.
\,In particular,~\,$\E{\bf U}=\E{\bf V}=0$.
\end{remark}

\begin{lemma} \label{2.2} {\rm  (Einmahl 1989, Lemma~5,~p.\ 55)}
Let \,$1\le m=(2\4s+1)\cdot2^r\le2^{N}$,
\,where \,$s,r$ are non-negative integers. Then
$$
S_m=\ffrac m{2^N}S_{2^N}+\sum_{n=r+1}^N\g_n\,\wt U_{n,l_{n,m}},
\eqno(2.41)
$$
where \,$\g_n=\g_n(m)\in[\40,\half\4]$ \,and the integers
~\,$l_{n,m}$ \,are defined by
$$
l_{n,m}\cdot2^n<m\le \big(l_{n,m}+1\big)\cdot2^n.
\eqno(2.42)
$$
\end{lemma}

The shortest proof of Lemma~\ref {2.2}
can be obtained with the help of a geometrical approach due
 to Massart~(1989,~p.\ 275).

\begin{remark}  \label{r2.2}
The inequalities\/~$(2.42)$ give a formal definition
of~\,$l_{n,m}$. \,To understand better
 the mechanism of the dyadic scheme,
one should remember another characterization of these numbers\/$:$
 \,$U_{n,l_{n,m}}$
\,is the sum over
the block of $2^n$ summands which contains~\,$X_m$,
\,the last summand in the sum~\,$S_m$.
\end{remark}

\begin{corollary} \label{c2.1}
Under the conditions of Lemma\/~{\rm \ref {2.2}}
$$
\bgl|S_m-T_m\bgr|\le
\bgl| U_{N,0}- V_{N,0}\bgr|
+\ffrac12\sum_{n=r+1}^N
\bgl|\wt U_{n,l_{n,m}}-\wt V_{n,l_{n,m}}\bgr|,\qquad
m=1,\dots,2^N.
$$
\end{corollary}

This statement evidently follows from Lemmas~\ref {l2.1}
and~\ref {2.2} and
from the relations~(2.9)--(2.12), (2.22) and~(2.23).

\vspace*{30pt}

\section {Proof of Theorem \ref{2.1}}\label{s3}

In the proof of Theorem~\ref {2.1} we shall use the following
auxiliary Lemmas
\ref{3.1}--\ref{3.4}
(Zaitsev 1995, 1996, 1998a).

\begin{lemma} \label{3.1} Suppose that \,$\L(\xi)\in\A$,
 \,$y\in{\bf R}^m$, \,$\a\in{\bf R}^1$.
\,Let \,${\bf M}:\Rd\to{\bf R}^m$ \,be a linear operator and
  \,$\wt\xi\in{\bf R}^k$ \,be the vector consisting of a subset
of coordinates of the vector~\,$\xi$. \,Then
\begin{eqnarray*}
\L({\bf M}\,\xi+y)\in{\cal A}_m\big(\norm {\bf M}\t\big), &&\qquad\hbox{where}
\quad\norm {\bf M}=\sup_{\norm x\le1}\norm{{\bf M}\4x},\\
\L(\a\4\xi)\in\AT{|\a|\4\t},&&\qquad\L(\wt\xi)\in\AD\t k.
\end{eqnarray*}
\end{lemma}

\begin{lemma} \label{3.2} Suppose that independent random vectors
\,$\xi^{(k)}$, \,$k=1,2$, \,satisfy the condition~
\,$\L\big(\xi^{(k)}\big)\in\AD\t{d_k}$.
\,Let \,$\xi=\bgl(\xi^{( 1)},\,\xi^{( 2)}\bgr)\in{\bf R}^{d_1+d_2}$
\,be the vector with
the first \,$d_1$~coordinates coinciding with those of~\,$\xi^{( 1)}$
and with the last \,$d_2$~coordinates
coinciding with those of~\,$\xi^{( 2)}$.
Then \,${\L(\xi)\in\AD\t{d_1+d_2}}$.
\end{lemma}

\begin{lemma} \label{3.3}
{\rm(Bernstein-type inequality)} \,Suppose that
 \,$\L(\xi)\in\AD\t1$, \,${\E\xi=0}$ \,and \,$\E\xi^2=\si^2$.
\,Then
$$
\P\bgl\{|\2\xi\2|
\ge x\bgr\}\le2\,\max\bgl\{\exp\bgl(-\pfrac{x^2\!}{4
\4\si^2}\bgr),
\,\exp\bgl(-\pfrac{x}{4\4\t}\bgr)\bgr\},\qquad x\ge0.
$$
\end{lemma}

\begin{lemma}  \label{3.4}
 Let the distribution of a random vector
~\,$\xi\in\Rd$ \,with \,$\E\xi=0$ \,satisfy the condition
\,$\L(\xi)\in\AV\t4d$, \,$\t\ge0$.
\,Assume that the variance \,${\si^2=\E\xi_d^2>0}$ \,of
the last coordinate $\xi_d$
of the vector $\xi$ is the minimal eigenvalue
of~\,$\cov \xi$.
\,Then there exist absolute positive constants
\,$c_8,\dots,c_{12}$
\,such that the following assertions hold$:$
\smallbreak

\noindent
{\rm a)} Let \,$d\ge2$.
\,Assume that \,$\xi_d$
\,is not correlated with previous coordinates \,$\xi_1,\dots,\xi_{d-1}$
\,of the vector~\,$\xi$.
\,Define
\,${\bf B}=\cov \,\ov{\bf P}_{d-1}\4\xi$ \,and denote by
\,$F(z\gm x)$,
\,$z\in{\bf R}^1$,
\, the r.c.d.f.\
 of \,$\xi_d$ \,for a given value
of~\,${\ov{\bf P}_{d-1}\4\xi=x\in{\bf R}^{d-1}}$.
\,Let
\,$\L(\ov{\bf P}_{d-1}\4\xi)\in\AV\t4{d-1}$.
\,Then there exists \,$y\in{\bf R}^1$
such that
$$
|y|\le c_8\,\t
\,\bgl\|{\bf B}\msqrt x\bgr\|^2
\le c_8\,\t\2\ffrac{\norm{x}^2}{\si^2},
\eqno(3.1)
$$
and
$$
\Phi_\si\big(\2z-\g(z)\2\big)<F(z+y\gm x)
<\Phi_\si\big(\2z+\g(z)\2\big),
\eqno(3.2)
$$
for
\,$\ftdsi\le c_9$,
\,${\bgl|{\bf B}\msqrt x\bgr|
\le \ffrac{c_{10}\4\si}{d^{3/2}\t}}$,
\,$|z|\le\ffrac{c_{11}\4\si^2}{d\4\t}$,
\,where $$
\g(z)=c_{12}\4\t\,\biggl(d^{3/2}
+d\4\de
\,\Bigl(1+\ffrac{|z|}{\si}\Bigr)
+\ffrac{z^2}{\si^2}\biggr),\qquad
\de=\bgl\|{\bf B}\msqrt x\bgr\|.
\eqno(3.3)
$$
\noindent
{\rm b)} The assertion\/ {\rm  a)} remains valid for
 \,$d=1$ \,with
\,$F(z\gm x)=\P\bgl\{\xi_1<z\bgr\}$
and \,$y=\de=0$ \,without any restrictions
on~\,${\bf B}$,~\,$\ov{\bf P}_{d-1}\4\xi$ \,and~\,$x$.
\end{lemma}

\begin{remark} \label{r3.1}
In\/ {\rm  Zaitsev} $(1995, 1996)$ the formulation of Lemma\/~{\rm \ref{3.4}}
is in some sense weaker,
see\/ {\rm  Zaitsev (1995, 1996,
Lemmas~$6.1$ and~$8.1)$}. In particular, instead of the conditions
$$
\L(\xi)\in\AV\t4d \qquad\hbox{and} \qquad
\L(\ov{\bf P}_{d-1}\4\xi)\in\AV\t4{d-1}
\eqno(3.4)
$$
the stronger conditions
$$
\L(\xi)\in\ATDRZ\t44\qquad\hbox{and} \qquad
\L(\ov{\bf P}_{d-1}\4\xi)\in\ADRZ\t44{d-1}
\eqno(3.5)
$$
are used.
However, in the proof of\/ $(3.1)$ and\/~$(3.2)$ only
the conditions\/~$(3.4)$
are applied.
The conditions\/ $(3.5)$ are necessary for the investigation
of quantiles of conditional distributions
corresponding to random vectors having coinciding
moments up to third order which has been
done in\/ {\rm  Zaitsev} $(1995, 1996)$
simultaneously with the proof of\/~$(3.1)$ and\/~$(3.2)$.
\end{remark}

\begin{lemma} \label{3.5}
Let \,$S_k=X_1+\dots+X_k$, \,$k=1,\dots,n$,
\,be sums of independent random vectors \,${X_j\in\Rd}$
\,and let \,$q(\cdt)$
\,be a semi-norm in~\,$\Rd$. Then
$$
\P\bgl\{\max_{1\le k\le n}\,q(S_k)>3\4t\bgr\}
\le3\,\max_{1\le k\le n}\P\bgl\{q(S_k)>t\bgr\},\qquad t\ge0.
\eqno(3.6)
$$
\end{lemma}

Lemma~\ref {3.5} is a version of the Ottaviani inequality,
see Dudley (1989, p.~251) or Hoffmann-J\o rgensen (1994,
p.~472). In the form (3.6) this inequality can be found in
Etemadi (1985) with 4 instead of~3 (twice).
The proof of Lemma~\ref {3.5} repeats those from the references
above and is therefore omitted.

\begin{lemma} \label{3.6}
Let the conditions of Theorem\/~{\rm \ref {2.1}} be satisfied
and assume that the vectors~\,$X_k$, \,$k=1,\dots,2^N$,
\,are constructed by the dyadic procedure
described in\/~{\rm (2.22)--(2.36).}
Then there exist
absolute positive constants~\,$c_{13},\dots,c_{17}$ \,such that
\smallbreak

\noindent{\rm a)}
If \,$\t\4d^{3/2}\!\big/2^{N/2}\le c_9$, then
$$
\bgl| U_{N,0}- V_{N,0}\bgr|\le
c_{13}\4d^{3/2}\4\t
\,\bgl(1+2^{-N}\4\bgl|U_{N,0}\bgr|^2\bgr)
\eqno(3.7)
$$
provided that
\,$\bgl|U_{N,0}\bgr|\le\ffrac{c_{14}\cdot2^N}{d^{3/2}\4\t};$
\smallbreak

\noindent{\rm b)}
If
\,$1\le n\le N$, \,$0\le k<2^{N-n}$,
\,$\t\4d^{3/2}\!\big/2^{n/2}\le c_{15}$,
\,then
$$
\bgl|\wt U_{n,k}-\wt V_{n,k}\bgr|\le
c_{16}\4d^{3/2}\4\t
\,\bgl(1+2^{-n}\4\bgl|{\bf U}_{n,k}\bgr|^2\bgr)
\eqno(3.8)
$$
provided that
\,$\bgl|{\bf U}_{n,k}\bgr|
\le\ffrac{c_{17}\cdot2^n}{d^{3/2}\4\t}.$
\end{lemma}

In the proof of Lemma~\ref {3.6} we need the following
auxiliary Lemma~\ref {3.7} which is useful for the application
of Lemma~\ref {3.4} to the conditional distributions
involved in the dyadic scheme.

\begin{lemma} \label{3.7} Let \,$F(\cdt)$ \,denote a continuous
distribution function and \,$G(\cdt)$ \,an arbitrary
distribution function
satisfying for \,$z\in B\in{\cal B}_1$ \,the inequality
$$
G\bgl(z-f(z)\big)<F(z+w)<G\bgl(z+f(z)\big)
$$
with some~\,$f:B\to{\bf R}^1$ \,and~\,$w\in{\bf R}^1$. \,Let~
 \,$\eta\in{\bf R}^1$, \,$0<G(\eta)<1$
\,and \,${\xi=F\me\big(G(\eta)\big)}$,
\,where \,${F\me(x)=
\sup\,\bgl\{u:F(u)\le x\bgr\}}$, \,$0<x<1$. \,Then
$$
|\4\xi-\eta\4|< f(\xi-w)+|\4w\4|,\qquad \hbox{if}
\quad \xi-w\in B.
$$
\end{lemma}
\noindent{\it Proof}\/ \ Put \,$\zeta=\xi-w$.
\,The continuity of~\,$F$
\,implies that \,$F\big(F\me(x)\big)\equiv x$, \,for \,${0<x<1}$.
\,Therefore,
$$
 \zeta\in B\Rightarrow
G\bgl(\zeta-f(\zeta)\big)<F(\xi)=G(\eta)\Rightarrow
\zeta-f(\zeta)<\eta\Rightarrow
\xi-\eta< f(\zeta)+w
$$
and
$$
 \zeta\in B\Rightarrow G(\eta)= F(\xi)<G\bgl(\zeta+f(\zeta)\big)
\Rightarrow\eta< \zeta+f(\zeta)\Rightarrow
\eta-\xi< f(\zeta)-w.
$$
This completes the proof of the lemma. $\square$\medbreak
\medskip

\noindent{\it Proof of Lemma\/ {\rm \ref{3.6}}} \ At first
we note that the conditions of Theorem~\ref {2.1}
imply that
$$
\cov {\bf U}_{n,k}
=2^n\4{\bf I}_{2d},\for1\le n\le N, \quad0\le k<2^{N-n},
$$
and, hence (see~(2.28)),
$$
\cov {\bf U}_{n,k}^j
=2^n\4{\bf I}_j,\for1\le j\le 2\4d.
\eqno(3.9)
$$

Let us prove the assertion~a). Introduce the vectors
$$
U_{N,0}^j=\bgl(U_{N,0}^{( 1)}\4,\dots,U_{N,0}^{(j)}\bgr),\qquad
V_{N,0}^j=\bgl(V_{N,0}^{( 1)}\4,\dots,V_{N,0}^{(j)}\bgr)
\eqno(3.10)
$$
consisting of the first $j$ coordinates of the vectors
$U_{N,0},\,V_{N,0}$ respectively. By (3.9), (2.32) and~(2.34),
$$
 U_{N,0}
=\ov{\bf P}_d\,{\bf U}_{N,0}
\eqno(3.11)
$$
and
$$
 U_{N,0}^j
={\bf U}_{N,0}^j,\qquad\cov U_{N,0}^j
=2^n\4{\bf I}_j,\for1\le j\le d.
\eqno(3.12)
$$
Moreover,
according to Lemma~\ref {l2.1}, Remark~\ref{r2.1},~(3.12) and~(2.19),
 the distributions $\L(U_{N,0}^j)$, \,$\jd$,
\,satisfy in the $j$-di\-men\-sional case the conditions
of Lemma~\ref {3.4} with \,$\si^2=2^N$ \,and
\,${\bf B}=\cov U_{N,0}^{j-1}=2^N\4{\bf I}_{j-1}$ (the last equality
for~\,${j\ge2}$).

Taking into account~(2.29) and
applying Lemmas~\ref{3.4} and~\ref{3.7},
 we obtain that
$$
\bgl|U_{N,0}^{( 1)}- V_{N,0}^{( 1)}\bgr|
\le c_{12}\4\t\,\Bigl(1+
\ffrac{\bgl|U_{N,0}^{( 1)}\bgr|^2}{2^N}\Bigr),
\eqno(3.13)
$$
if \,$\ffrac\t{2^{N/2}}\le c_9$,
\,$\bgl|U_{N,0}^{( 1)}\bgr|\le\ffrac{c_{11}\cdot2^N}{\t}$.
\,Furthermore,
$$
\bgl|U_{N,0}^{(j)}- V_{N,0}^{(j)}\bgr|
\le c_{12}\4\t\,\biggl(j^{3/2}
+j^{3/2}\4\ffrac{\bgl|U_{N,0}^{j-1}\bgr|}{2^{N/2}}
\,\Bigl(1+\ffrac{\bgl|U_{N,0}^{(j)}-y_j\bgr|}{2^{N/2}}\Bigr)
$$
$$\kern6cm
+\ffrac{\bgl|U_{N,0}^{(j)}-y_j\bgr|^2}{2^N}
\biggr)+|\4y_j\4|,
\eqno(3.14)
$$
if
$$
\ffrac{\t\4j^{3/2}}{2^{N/2}}\le c_9, \quad
\ffrac{\bgl|U_{N,0}^{j-1}\bgr|}{2^{N/2}}
\le \ffrac{c_{10}\cdot2^{N/2}}{j^{3/2}\t},
\quad
\bgl|U_{N,0}^{(j)}-y_j\bgr|\le\ffrac{c_{11}\cdot2^N}{j\4\t}, \qquad
2\le j\le d,
\eqno(3.15)
$$
where
$$
|\4y_j\4|\le c_8\4\t\4j\4
\ffrac{\bgl|U_{N,0}^{j-1}\bgr|^2}{2^N},
\qquad2\le j\le d.
\eqno(3.16)
$$
Obviously,
$$
\bgl|U_{N,0}^{( 1)}\bgr|\le
\max\bgl\{\bgl|U_{N,0}^{j-1}\bgr|,\,\bgl|U_{N,0}^{(j)}\bgr|\bgr\}
=\bgl|U_{N,0}^j\bgr|\le\bgl|U_{N,0}\bgr|, \qquad 2\le j\le d,
\eqno(3.17)
$$
see~(2.31) and (3.10).
Using (3.13), (3.14), (3.16) and (3.17), we see that one can choose
 $c_{13}$ to be so large and $c_{14}$ to be so small that
$$
\bgl|U_{N,0}^{(j)}- V_{N,0}^{(j)}\bgr|
\le c_{13}\4d^{3/2}\4\t
\,\bgl(1+2^{-N}\4\bgl|U_{N,0}\bgr|^2\bgr),
\eqno(3.18)
$$
if
\,$
\ffrac{\t\4d^{3/2}}{2^{N/2}}\le c_9$, \,$
\bgl|U_{N,0}\bgr|
\le \ffrac{c_{14}\cdot2^N}{d^{3/2}\t}$,
\,$1\le j\le d$.
\,The inequality (3.7) immediately follows
from~(3.18),~(2.23) and~(2.31).

Now we shall prove item b).
According to Lemma~\ref {l2.1}, Remark~\ref{r2.1}, (2.18), (2.31) and~(3.9),
the distributions $\L({\bf U}_{n,k}^j)$, \,$j=d+1,\dots,2\4d$,
\,satisfy in the $j$-di\-men\-sional case the conditions
of Lemma~\ref {3.4} with \,$\si^2=2^n$,
 \,${{\bf B}=\cov {\bf U}_{n,k}^{j-1}=2^n\4{\bf I}_{j-1}}$.

Using~(2.33) and
applying Lemmas~\ref{3.4} and~\ref{3.7}, we obtain that
$$
\bgl|{\bf U}_{n,k}^{(j)}-{\bf V}_{n,k}^{(j)}\bgr|
\le c_{12}\4\t\,\biggl(j^{3/2}
+j^{3/2}\4\ffrac{\bgl|{\bf U}_{n,k}^{j-1}\bgr|}{2^{n/2}}
\,\Bigl(1+\ffrac{\bgl|{\bf U}_{n,k}^{(j)}-y_j\bgr|}{2^{n/2}}\Bigr)
$$
$$
\kern6cm
+\ffrac{\bgl|{\bf U}_{n,k}^{(j)}-y_j\bgr|^2}{2^n}
\biggr)+|\4y_j\4|,
\eqno(3.19)
$$
if
$$
\ffrac{\t\4j^{3/2}}{2^{n/2}}\le c_9, \quad
\ffrac{\bgl|{\bf U}_{n,k}^{j-1}\bgr|}{2^{n/2}}
\le \ffrac{c_{10}\cdot2^{n/2}}{j^{3/2}\t},
\quad
\bgl|{\bf U}_{n,k}^{(j)}-y_j\bgr|\le\ffrac{c_{11}\cdot2^n}{j\4\t},
\eqno(3.20)
$$
where
$$
|\4y_j\4|\le c_8\4\t\4j\4
\ffrac{\bgl|{\bf U}_{n,k}^{j-1}\bgr|^2}{2^n},
\qquad d+1\le j\le 2\4d.
\eqno(3.21)
$$
Obviously,
$$
\max\bgl\{\bgl|{\bf U}_{n,k}^{j-1}\bgr|,
\,\bgl|{\bf U}_{n,k}^{(j)}\bgr|\bgr\}
=\bgl|{\bf U}_{n,k}^j\bgr|\le\bgl|{\bf U}_{n,k}\bgr|,
\eqno(3.22)
$$
see~(2.34).
Using (3.19), (3.21) and (3.22), we see that one can choose
 $c_{15}$ and $c_{17}$ to be so small
and $c_{16}$ to be so large that
$$
\bgl|{\bf U}_{n,k}^{(j)}-{\bf V}_{n,k}^{(j)}\bgr|
\le c_{16}\4d^{3/2}\4\t
\,\bgl(1+2^{-n}\4\bgl|{\bf U}_{n,k}\bgr|^2\bgr)\nopagebreak
\eqno(3.23)
$$
if \,$ \ffrac{\t\4d^{3/2}}{2^{n/2}}\le c_{15}$, \,$ \bgl|{\bf
U}_{n,k}\bgr| \le \ffrac{c_{17}\cdot2^n}{d^{3/2}\t}$, \,$ d+1\le
j\le 2\4d$. \,The inequality (3.8) immediately follows
from~(3.23),~(2.24),~(2.25) and~(2.34). $\square$\medbreak

{\it Proof of Theorem\/~{\rm \ref {2.1}}} \
Let~\,$X_k$, \,$k=1,\dots,2^N$, \,denote the vectors
 constructed by the dyadic procedure
described in~(2.22)--(2.36).
Denote
$$
\DE=\DE(2^N)=\max_{1\le k\le2^N}\,\bgl|S_k-T_k\bgr|,
\eqno(3.24)
$$
$$
 c_5
=\min\,\bgl\{c_9,\,c_{15}\bgr\},
\quad
 c_{18}
=\min\,\bgl\{c_{14},\,c_{17},\,1\bgr\},
\quad
y\=\ffrac{c_{18}}{d^{3/2}\4\t}\le\ffrac1\t,
\eqno(3.25)
$$
fix some \,$x>0$ \,and choose the integer~\,$M$ \,such that
$$
x<4\4y\cdot2^M\le2\4x.
\eqno(3.26)
$$

We shall estimate \,$\P\bgl\{\DE\ge x\bgr\}$.
\,Consider separately two possible cases: \,$M\ge N$
\,and \,$M< N$. \,Let, at first, \,$M\ge N$. \,Denote
$$
\DE_1=\max_{1\le k\le2^N}\,\bgl|S_k\bgr|,\qquad
\DE_2=\max_{1\le k\le2^N}\,\bgl|T_k\bgr|.
\eqno(3.27)
$$
It is easy to see that
\,$
\DE\le\DE_1+\DE_2
$ \,and, hence,
$$
\P\bgl\{\DE\ge x\bgr\}\le\P\bgl\{\DE_1\ge x/2\bgr\}
+\P\bgl\{\DE_2\ge x/2\bgr\}.
\eqno(3.28)
$$
Taking into account the completeness of classes \,$\A$
\,with respect to convolution,
 applying Lemmas~\ref{3.5}, \ref{3.1} and ~\ref{3.3}
and using~(3.25) and~(3.26), we obtain
that \,$2^N\le 2^M\le x/2\4y$ \,and
\begin{eqnarray*}\vspace{5pt}
 \P\bgl\{\DE_1\ge x/2\bgr\}\5
&\le& \5 3\,\max\limits_{1\le k\le2^N}\,\P\bgl\{\bgl|S_k\bgr|\ge x/6\bgr\}\5\\
&\le& \5 6\4d\,\exp\Big(-\min\Big\{\ffrac{x^2}{144\cdot2^N},
\ffrac{x}{24\4\t}\Big\}\Big)\\
&\le& \5 6\4d
\,\exp\Big(-\ffrac{c_{19}\,x}{d^{3/2}\4\t}\Big).
\hbox{\rlap{\hskip3.9cm(3.29)}}
\end{eqnarray*}
Since all $d$-dimensional Gaussian distributions
 belong to all classes~\,$\A$, \,${\t\ge0}$,
\,we~automatically obtain that
$$
\P\bgl\{\DE_2\ge x/2\bgr\}
\le6\4d
\,\exp\Big(-\ffrac{c_{19}\,x}{d^{3/2}\4\t}\Big).
\eqno(3.30)
$$
From (3.28)--(3.30) it follows in the case \,$M\ge N$ \,that
$$
\P\bgl\{\DE\ge x\bgr\}
\le 12\4d\,\exp\Big(-\ffrac{c_{19}\,x}{d^{3/2}\4\t}\Big).
\eqno(3.31)
$$

Let now \,$M< N$. \,Denote
$$
L=\max\bgl\{0,\,M\bgr\}
\eqno(3.32)
$$
and
$$
\begin{tabular}{r c l r}\vspace{5pt}
$\DE_3\5$&=&$\5\max\limits_{0\le k<2^{N-L}}\,\max\limits_{1\le l\le2^L}
\,\bgl|S_{k\cdot2^L+l}-S_{k\cdot2^L}\bgr|,$&\hfill
\rlap{\kern2.2cm
(3.33)}\\ \vspace{5pt}
$\DE_4\5$&=&$\5\max\limits_{0\le k<2^{N-L}}\,\max\limits_{1\le l\le2^L}
\,\bgl|T_{k\cdot2^L+l}-T_{k\cdot2^L}\bgr|,$ &\hfill
\rlap{\kern2.2cm
(3.34)}\\
$\DE_5\5$&=&
$\5\max\limits_{1\le k\le2^{N-L}}
\,\bgl|S_{k\cdot2^L}-T_{k\cdot2^L}\bgr|.$&\hfill
\rlap{\kern2.2cm
(3.35)}
\end{tabular}
$$
Introduce the event
$$
A=\bgl\{\4\om:\,\bgl|U_{L,k}\bbgr|
<y\cdot2^L, \
0\le k<2^{N-L}\4\bgr\}
\eqno(3.36)
$$
(we assume that all considered random vectors
are measurable mappings of~\,${\om\in\Omega}$).
For the complementary event we use the notation
~\,$\ovln A=\Omega\setminus A$.

We consider separately two possible cases: \,$L=M$
\,and \,$L=0$. \,Let~\,${L=M}$.
\,It is evident that in this case
$$
\DE\le\DE_3+\DE_4+\DE_5.
\eqno(3.37)
$$
 Moreover, by virtue of (3.37), (3.26), (3.33) and~(3.36),
we have
$$
\ovln A\subset\bgl\{\om:\DE_3\ge x/4\bgr\}.
\eqno(3.38)
$$
From (3.37) and (3.38) it follows that
$$
{}\P\bgl\{\DE\ge x\bgr\}\le\P\bgl\{\DE_3\ge x/4\bgr\}
+\P\bgl\{\DE_4\ge x/4\bgr\}+\P\bgl\{\DE_5\ge x/2, \,A\bgr\}.
\eqno(3.39)
$$
Using Lemmas~\ref{3.5}, \ref{3.1} and ~\ref{3.3},
the completeness of classes \,$\A$
\,with respect to convolution
 and  the relations
~(3.25) and~(3.26), we obtain, for \,$0\le k<2^{N-L}$,
that~\,${2^L= 2^M\le x/2\4y}$ \,and
\begin{eqnarray*}
{}\P\bgl\{\max_{1\le l\le2^L}
\,\bgl|S_{k\cdot2^L+l}-S_{k\cdot2^L}\bgr|\ge x/4\bgr\}\5
&\le&\5 3\,\max_{1\le l\le2^L}\,\P\bgl\{
\bgl|S_{k\cdot2^L+l}-S_{k\cdot2^L}\bgr|\ge x/12\bgr\}
\\
&\le&\5 6\4d\,\exp\Big(-\min\Big\{\ffrac{x^2}{576\cdot2^L},
\ffrac{x}{48\4\t}\Big\}\Big)
\\
&\le&\5{6\4d
\,\exp\Big(-\ffrac{c_{20}\4x}{d^{3/2}\4\t}\Big).
\hskip2.8cm(3.40)}
\end{eqnarray*}
Since all $d$-dimensional Gaussian distributions
 belong to classes~\,$\A$ \,for all \,$\t\ge0$,
\,we immediately obtain that
$$
{}\P\bgl\{\max_{1\le l\le2^L}
\,\bgl|\4T_{k\cdot2^L+l}-T_{k\cdot2^L}\bgr|\ge x/4\bgr\}
\le6\4d
\,\exp\Big(-\ffrac{c_{20}\4x}{d^{3/2}\4\t}\Big).
\eqno(3.41)
$$
From (3.33), (3.34), (3.40) and (3.41) it follows that
$$
{}\P\bgl\{\DE_3\ge x/4\bgr\}
+\P\bgl\{\DE_4\ge x/4\bgr\}
\le2^N\cdot12\4d
\,\exp\Big(-\ffrac{c_{20}\4x}{d^{3/2}\4\t}\Big).
\eqno(3.42)
$$

Assume that \,$L=0$. \,Then, according to (3.24) and~(3.35),
 \,$\DE=\DE_5$ \,and, hence, we have the rough bound
$$
{}\P\bgl\{\DE\ge x\bgr\}\le\P\bgl\{\ovln A\bgr\}
+\P\bgl\{\DE_5\ge x/2, \,A\bgr\}.
\eqno(3.43)
$$
In this case \,$U_{L,k}=X_{k+1}$, \,$2^L=1\ge2^M$,
\,$y> x/4$ \,(see~(3.25), (3.26) and~(3.32)).
Therefore, by~(3.36) and by~Lemmas~\ref{3.1} and~\ref{3.3},
\begin{eqnarray*}
{}\P\bgl\{\ovln A\bgr\}&\le&\sum\limits_{k=0}^{2^N-1}
\P\bgl\{\bgl|U_{L,k}\bbgr|
\ge y\cdot2^L\bgr\}=
\sum\limits_{k=1}^{2^N}
\P\bgl\{\bgl|X_{k}\bbgr|
\ge{y}\bgr\}
\\
&\le&2^{N+1}\4d\,\exp\Big(-\min\Big\{\ffrac{y^2}{4},
\ffrac{y}{4\,\t}\Big\}\Big)\\
&\le&2^{N+1}\4d\,\exp\Big(-\min\Big\{\ffrac{x\4y}{16},
\ffrac{x}{16\,\t}\Big\}\Big)\\
&\le&2^{N+1}\4d
\,\exp\Big(-\ffrac{c_{21}\4x}{d^{3/2}\4\t}\Big).
\hbox{\rlap{\hskip4.37cm(3.44)}}
\end{eqnarray*}

It remains to estimate \,$\P\bgl\{\DE_5\ge x/2, \,A\bgr\}$
\,in both cases: \,$L=M$ \,and \,${L=0}$
\,(see~(3.39) and (3.42)--(3.44)).
Let \,$L$ \,defined
by~(3.32) be arbitrary.
Fix an integer~\,$k$ \,satisfying \,${1\le k\le2^{N-L}}$
\,and denote for simplicity
$$
j=j(k)\=k\cdot2^{L}.
\eqno(3.45)
$$
By Corollary~\ref{c2.1},
\,we have
$$
\bgl|S_{k\cdot2^{L}}-T_{k\cdot2^{L}}\bbgr|=
\bgl|S_{j}-T_{j}\bgr|\le
\bgl| U_{N,0}- V_{N,0}\bgr|
+\ffrac12\sum_{n=L+1}^N
\bgl|\wt U_{n,l_{n,j}}-\wt V_{n,l_{n,j}}\bgr|,
\eqno(3.46)
$$
where \,$l_{n,j}$ \,are integers, defined by
\,$l_{n,j}\cdot2^n<j\le \big(l_{n,j}+1\big)\cdot2^n\,$
(see~(2.42)).

By virtue of (3.25) and~(3.36), for \,$\om\in A$ \,we have
$$
\bgl|U_{L,l}\bbgr|
<y\cdot2^L
=\ffrac{c_{18}\cdot2^L}
{d^{3/2}\4\t}\le\ffrac{\min\{c_{14},\,c_{17}\}\cdot2^L}
{d^{3/2}\4\t},\qquad 0\le l<2^{N-L},
\eqno(3.47)
$$
and, by~(2.35)--(3.37),
 \,$U_{L,l}$ \,are sums over
blocks consisting of \,$2^L$ \,summands. Moreover, \,$ U_{n,l}$
\,(resp. \,$\wt U_{n,l}$),
\,$L+1\le n\le{N}$, \,$0\le l<2^{N-n}$,
 \,are sums (resp. differences)
of two sums over blocks containing each $2^{n-1}$ summands.
These sums and differences can be represented as linear combinations
(with coefficients~\,$\pm1$) of \,$2^{n-L}$
\,sums over blocks
containing each \,$2^L$~\,summands and satisfying~(3.47).
Therefore, for \,$\om\in A$, \,$L+1\le n\le{N}$,
\,$0\le l<2^{N-n}$ \,we have
(see~(2.32) and~(2.34))
$$
\bgl|{\bf U}_{n,l}\bbgr|=
\max\bgl\{\bgl|U_{n,l}\bbgr|, \,\bgl|\wt U_{n,l}\bbgr|\bgr\}
\le2^{n-L}\4y\cdot2^L
=y\cdot2^n\le\ffrac{\min\{c_{14},\,c_{17}\}\cdot2^n}
{d^{3/2}\4\t}.
\eqno(3.48)
$$
Using (3.48), we see that if \,$\om\in A$,
\,the conditions of Lemma~\ref {3.6} are satisfied
for \,$\t$, \,$ U_{N,0}$ \,and~
\, ${\bf U}_{n,l}$, \,if \,$L+1\le n\le{N}$, \,$0\le l<2^{N-n}$.
By (3.46), (3.48) and by Lemma~\ref {3.6},
for~\,$\om\in A$ \,we have
$$
\bgl|S_{j}-T_{j}\bbgr|\le
c_{13}\4d^{3/2}\4\t
\,\bgl(1+2^{-N}\4\bgl|U_{N,0}\bgr|^2\bgr)\hskip6cm{}
$$
$$
+\sum_{n=L+1}^N
c_{16}\4d^{3/2}\4\t
\,\Big(\,1+2^{-n}\4\max\bgl\{\bgl|U_{n,l_{n,j}}\bgr|^2,
\,\bgl|\wt U_{n,l_{n,j}}\bgr|^2\bgr\}\Big)
$$
$$
\le
c\4d^{3/2}\4\t\,\biggl(
\,N+1+2^{-N}\4\bgl|U_{N,0}\bgr|^2
+\sum_{n=L}^{N-1}
2^{-n}\4\bgl(\bgl|\YY U{}n\bgr|^2+
\bgl|U_{(n)}\bgr|^2\bgr)\biggr),
\eqno(3.49)
$$
where
$$
\YY U{}n=U_{n,l_{n,j}},\qquad
 U_{(n)}=U_{n,\wt l_{n,j}},
\eqno(3.50)
$$
and
$$
\wt l_{n-1,j}=
\left\{
\begin{tabular}{c l}
$2\4l_{n,j},$&\quad$\hbox{if}\quad l_{n-1,j}=2\4l_{n,j}+1,$\\
$2\4l_{n,j}+1,$&\quad$\hbox{if}\quad l_{n-1,j}=2\4l_{n,j},$
\end{tabular}
\right.
\qquad L< n\le{N}
\eqno(3.51)
$$
(it is easy to see that \,$l_{n-1,j}$ \,can be equal either to
\,$2\4l_{n,j}$ \,or to \,$2\4l_{n,j}+1$, \,for given~\,$l_{n,j}$).
In other words, \,$\YY U{}n$, \,$L\le n\le{N}$, \,is the sum over
the block of~$2^n$~summands which contains~\,$X_{j}$. \,The
sum~\,$U_{(n)}$ \,does not contain~\,$X_{j}$ \,and
$$
\YY U{}{n+1}=\YY U{}n+U_{(n)},
\qquad L\le n<{N}\nopagebreak
\eqno(3.52)
$$
(see~(3.37)). The equality (3.52) implies
$$
\YY U{}n=
 \YY U{}{L}+
\sum_{s=0}^{n-L-1}U_{(L+s)},\qquad L\le n\le{N}.
\eqno(3.53)
$$
It is important that all summands in the right-hand side
of~(3.53) are the sums of disjoint blocks of independent summands.
Therefore, they are independent.

Put \,$\be=1/\sqrt2$. \,Then, using (3.53) and
 the H\"older inequality, one can easily derive
 that,
for \,$ L\le n\le{N}$,
$$
\bgl|\YY U{}n\bbgr|^2\le
c_{22}\,\bigg(\,\be^{-(n-L)}\4 \bgl|\YY U{}{L}\bbgr|^2+
\sum_{s=0}^{n-L-1}\be^{-(n-L-1)+s}\4
\bgl|U_{(L+s)}\bbgr|^2\,\bigg),
\eqno(3.54)
$$
with
\,$c_{22}=\sum\limits_{j=0}^{\infty}\be^{j}=\ffrac{\sqrt2}{\sqrt2-1}$.
\,It is easy to see that
$$
\sum_{n=L}^{N}2^{-n}
\,\be^{-(n-L)}\4 \bgl|\YY U{}{L}\bbgr|^2\le
c_{22}\cdot2^{-L}\4 \bgl|\YY U{}{L}\bbgr|^2.
\eqno(3.55)
$$
Moreover,
$$
\sum_{n=L+1}^{N}\sum_{s=0}^{n-L-1}
2^{-n}\4\be^{-(n-L-1)+s}\4
\bgl|U_{(L+s)}\bbgr|^2\hskip3cm{}
$$
$$
=\sum_{s=0}^{N-L-1}\sum_{n=L+1+s}^{N}
2^{-n}\4\be^{-(n-L-1)+s}\4
\bgl|U_{(L+s)}\bbgr|^2$$
$$
{}\hskip3cm\le c_{22}\sum_{s=0}^{N-L-1}2^{-(L+1+s)}\4
\bgl|U_{(L+s)}\bbgr|^2.\qquad
\eqno(3.56)
$$
It is clear that the inequalities (3.54)--(3.56) imply
$$
2^{-N}\4\bgl|U_{N,0}\bgr|^2+\sum_{n=L}^{N-1}
2^{-n}\4\bgl(\bgl|\YY U{}n\bgr|^2+
\bgl|U_{(n)}\bgr|^2\bgr)\hskip3cm{}
$$
$$
\le c_{22}\,\bigg(\ffrac{ \bgl|\YY U{}{L}\bbgr|^2}{2^{L}}
+\sum_{s=0}^{N-L-1}\ffrac{
\bgl|U_{(L+s)}\bbgr|^2}{2^{L+1+s}}\bigg)
+\sum_{n=L}^{N-1}\ffrac{\bgl|U_{(n)}\bgr|^2}{2^{n}}
$$
$$
\hskip3cm\le c\,\bigg(\ffrac{ \bgl|\YY U{}{L}\bbgr|^2}{2^{L}}
+\sum_{n=L}^{N-1}\ffrac{\bgl|U_{(n)}\bgr|^2}{2^{n}}\bigg).
\eqno(3.57)
$$
\smallskip

\noindent
From (3.49) and (3.57) it follows that
for \,$\om\in A$ \,we have
$$
\bgl|S_{j}-T_{j}\bbgr|\le
c_{23}\4d^{3/2}\4\t\,\biggl(\,N+1+
\ffrac{ \bgl|\YY U{}{L}\bbgr|^2}{2^{L}}
+\sum_{n=L}^{N-1}\ffrac{\bgl|U_{(n)}\bgr|^2}{2^{n}}\bigg).
\eqno(3.58)
$$

Denote (for \,$0\le n\le N$, \,$0\le l<2^{N-n}$)
$$
W_{n,l}=\left\{\begin{tabular}{c l}
$2^{-n}\4\bgl|U_{n,l}\bgr|^2,$\quad&\hbox{if}
\quad$\bgl|U_{n,l}\bgr|\le y\cdot2^n,$\\
 \qquad 0,&\hbox{otherwise}.\end{tabular}\right.
\eqno(3.59)
$$
Let us show that
$$
\E\exp\big(t\,W_{n,l}\big)\le2\4d+1,\for0\le t\le \ffrac1{8}.
\eqno(3.60)
$$
Indeed, integrating by parts, we obtain
$$
\E\exp(t\4W_{n,l})=1+\int\limits_0^{y^2\cdot2^n}
\,t\4\exp(t\4u)\,\P\bgl\{W_{n,l}\ge u\bgr\}\,du\hskip2.45cm{}
$$
$$
{}\hskip2.45cm\le1+\ffrac18\int\limits_0^{y^2\cdot2^n}
\,\4\exp\big(u/8\big)\,\P\bgl\{\bgl|U_{n,l}\bgr|
\ge 2^{n/2}\4\sqrt u\bgr\}\,du.
\eqno(3.61)
$$
Taking into account (3.37), (3.25) and
using Lemmas~\ref{3.1} and~\ref{3.3}, we obtain that
$$
{\P\bgl\{\bgl|U_{n,l}\bgr|}
\ge 2^{n/2}\4\sqrt u\bgr\}
\le2\4d\,\exp\Big(-\min\Big\{\ffrac{2^n\4u}{4\cdot2^n},
\ffrac{2^{n/2}\4\sqrt u}{4\,\t}\Big\}\Big)
$$
$$
{}\hskip2.37cm{}
\le2\4d\,\exp\Big(-\min\Big\{\ffrac{u}{4},
\ffrac{u}{4\4y\4\t}\Big\}\Big)
$$
$${}\hskip.24cm{}
=2\4d
\,\exp\Big(-\ffrac{u}{4}\Big),
\eqno(3.62)
$$
if \,$0\le u\le y^2\cdot2^n$.
\,The relation~(3.60) immediately follows from
~(3.61) and~(3.62).

The relations (3.47), (3.48) and~(3.59) imply that,
for \,$ L\le n\le{N}$, \,$0\le l<2^{N-n}$,
\,$\om\in A$,
$$
2^{-n}\4\bgl|U_{n,l}\bgr|^2=W_{n,l}.
\eqno(3.63)
$$
Thus, according to (3.50), we can rewrite (3.58) in the form
$$
\bgl|S_{j}-T_{j}\bgr|\le
c_{23}\4d^{3/2}\4\t\,\biggl(\,N+1+\YY W{}L
+\sum_{n=L}^{N-1}W_{(n)}\bigg),\qquad \om\in A,\nopagebreak
\eqno(3.64)
$$
where
$$
\YY W{}L=W_{L,l_{L,j}},\qquad
 W_{(n)}=W_{n,\wt l_{n,j}},
\eqno(3.65)
$$

Putting now
\,$t^*=(8\,c_{23}\4d^{3/2}\4\t)\me$ \,and
\,${t=t^*\cdot c_{23}\4d^{3/2}\4\t=1/8}$,
\,taking into account that the random variables
\,$\YY W{}L$, $W_{(L)}$, \dots,
 $W_{(N-1)}$
\,are independent and applying
~(3.60),~(3.64) and~(3.65), we obtain
\begin{eqnarray*}
\P\Bigl\{\bgl\{&&\hskip-.7cm\om:
\bgl|S_{j}-T_{j}\bbgr|
\ge x/2\bgr\}\cap A\Bigr\}\\
&&\le\P\Bigl\{
\,c_{23}\4d^{3/2}\4\t\,\Bigl(\,N+1+\YY W{}L
+{\sum\limits_{n=L}^{N-1}}W_{(n)}\,\Big)
\ge x/2\,\Bigr\}\\
&&\le\P\Bigl\{t\,\Bigl(\,\YY W{}L
+{\sum\limits_{n=L}^{N-1}}
W_{(n)}\,\Big)\ge t^*x/2-t\4(N+1)\Bigr\}\\
&&\le\E\exp\biggl(t\,\Bigl(\,\YY W{}L
+{\sum\limits_{n=L}^{N-1}}W_{(n)}\,\Big)\bigg)
\Big/\exp\bgl( t^*x/2-t\4(N+1)\bgr)
\\
&&=\E\exp\bgl(t\4\YY W{}L\bgr)
{\prod\limits_{n=L}^{N-1}}
\E\exp\bgl(t\4W_{(n)}\bgr)
\Big/\exp\bgl( t^*x/2-t\4(N+1)\bgr)\\
&&\le\,(3\4d)^{N+1}\,\exp\Big(\ffrac {N+1}{8}
-\ffrac x{16\,c_{23}\4d^{3/2}\4\t}\Big).
\hskip4.05cm(3.66)
\end{eqnarray*}
From (3.35), (3.45) and (3.66) it follows that
$$
{}\P\bgl\{\DE_5\ge x/2, \,A\bgr\}
\le2^N\cdot(3\4d)^{N+1}\,\exp\Big(\ffrac {N+1}{8}
-\ffrac x{16\,c_{23}\4d^{3/2}\4\t}\Big).
\eqno(3.67)
$$
Using (3.31), (3.39), (3.42)--(3.44) and (3.67), we obtain that
$$
{}\P\bgl\{\DE\ge x\bgr\}\le(19\4d)^{N+1}
\,\exp\Big(-\ffrac x{c_{24}\4d^{3/2}\4\t}\Big),
\qquad x\ge0,
\eqno(3.68)
$$
where we can take
\,$
c_{24}=\max
\,\bgl\{16\4c_{23},\, c_{19}\me,\, c_{20}\me,\, c_{21}\me,\, 2\bgr\}
$.
\,Let the quantities \,$\e, \,x_0>0$ \,be defined by the relations
$$
\e=\ffrac 1{2\4c_{24}\4d^{3/2}\4\t}\le
\ffrac 1{4\4\t},
\qquad
 e^{\e x_0}=(19\4d)^{N+1}.
\eqno(3.69)
$$
Integrating by parts and using (3.68) and (3.69), we obtain
$$
\E e^{\e\DE}=\int _0^\infty \e\4e^{\e x}
\,\P\bgl\{\DE\ge x\bgr\}\,dx+1,
$$
$$
\int _0^{x_0} \e\4e^{\e x}
\,\P\bgl\{\DE\ge x\bgr\}\,dx\le
\int _0^{x_0} \e\4e^{\e x}\,dx=e^{\e x_0}-1
=(19\4d)^{N+1}-1,
$$
$$
\int _{x_0}^\infty \e\4e^{\e x}
\,\P\bgl\{\DE\ge x\bgr\}\,dx\le
\int _{x_0}^\infty \e\4e^{-\e (x-x_0)}\,dx=1,
$$
and, hence,
$$
\E e^{\e\DE}\le (19\4d)^{N+1}+1\le (20\4d)^{N+1}.
$$
Together with (3.24) and (3.69), this completes the proof of
Theorem~\ref {2.1}. $\square$\medbreak

\vspace*{30pt}

\section{Proofs of Theorems \ref {1.1}--\ref{1.4}}\label{s4}

We start the proofs of Theorems \ref {1.1}--\ref {1.3}  with
the following common part.
\medbreak

\noindent{\it Beginning of the proofs of Theorems~
{\rm\ref{1.1}},\/ {\rm \ref{1.2}} and\/ {\rm \ref{1.3}}} \
At first we shall verify that under the conditions of
 Theorems~\ref {1.2}  or~\ref {1.3}  we have
\,${\L(\xi_k)\in \A}$.
\,For Theorem~\ref {1.3}  this relation
is an immediate consequence of Lemma~\ref {3.1}, of
the completeness of classes~\,$\A$ \,with respect to convolution
and of the conditions~(1.8) and~(1.10)--(1.12).
In the case of Theorem~\ref {1.2}  we denote \,$K=\L(\eta)$.
\,One can easily verify that \,${\bf B}=\cov K=\gamma^2\,{\bf I}_d$,
\,where \,$\gamma^2$ \,is defined by~(1.7) and, hence,
$$
1\le\gamma^2\le3.
\eqno(4.1)
$$
Moreover,
$$
\p(K,z)=\log\E e^{\8z,\eta\9}=
\log\ffrac{\bgl(4+\t^2\4(d+\8z,\ov z\9)\bgr)
\,\exp\big(\8z,\ov z\9\!/2\big)}
{(4+\t^2\4d)},\qquad z\in\Cd.
\eqno(4.2)
$$
Using~(4.1) and~(4.2), we obtain
$$
\bgl|d_u d_v^2\4\p(K,z)\bgr|
=\bgl|d_u d_v^2\,\log \big(4+\t^2\4(d+\8z,\ov z\9)\big)\bgr|\le
c\4\t^3\nnnorm u\nnnorm v^2\le\nnnorm u\t\<{\bf B} \,v,v\>,
\eqno(4.3)
$$
for \,$\nnnorm z\t\le1$, \,provided that \,$c_1$ \,(involved in
Assertion~A) is sufficiently small. \,This means
that~\,${K=\L(\eta)\in\A}$. The relation
\,$\L(\xi_k)=\L\big(\eta/\gamma\big)\in\A$, \,${k=1,\dots,n}$,
\,follows from~(4.1) and from Lemma~\ref {3.1}.

The text below is related to
Theorems \ref {1.1}, \ref {1.2}  and \ref {1.3}
simultaneously.
Without loss of generality we assume that the amount
of summands is equal to \,$2^N$
with some positive integer~\,$N$.
\,It suffices to show that the dyadic scheme related to the vectors
\,$\xi_1,\dots,\xi_{2^N}$ \,satisfies the conditions of Theorem~\ref {2.1}
with \,$\t^*=\sqrt2\,\t$ \,instead of~\,$\t$. \,According to Lemma~\ref {l2.1},
we can verify the conditions~(2.18) and~(2.19)
for the vectors~\,${\bf U}_{n,k}^j$ and~\,${\bf U}_{N,0}^j$
\,instead of~\,${\bf U}_{n,k}^{*j}$ and~\,${\bf U}_{N,0}^{*j}$.
\,To this end we shall show that
$$
{\cal L}\big({\bf U}_{n,k}^j\big)\in \AV{\sqrt 2\,\t}4j
\for 0\le k<2^{N-n},\quad1\le n\le N,\quad
1\le j\le 2\4d.
\eqno(4.4)
$$

Recall that \,${\bf U}_{n,k}={\bf A}\4\wt{\bf U}_{n,k}$,
\,where \,${\bf A}$ \,is the linear operator defined
by~(2.16) and satisfying~(2.40).
Furthermore, \,$\wt{\bf U}_{n,k}
=\bgl(U_{n-1,2k},\,U_{n-1,2k+1})\in{\bf R}^{2d}$,
\,where the $d$-dimensional vectors
\,$U_{n-1,2k}$ \,and~\,$U_{n-1,2k+1}$
\,are independent.
The relation
\,${\L({\bf U}_{n,k})\in\AD{\sqrt 2\,\t}{2d}}$
\,can be therefore easily derived from
the conditions of Theorems \ref {1.1},
\ref {1.2}  and~\ref {1.3}  with the help of
Lemmas~\ref {l2.1},~\ref{3.1} and~\ref{3.2} (see~(2.40))
if we take into account the completeness
of classes~\,$\A$ with respect to convolution
and their monotonicity with respect to~\,$\t$.
It is easy to see
 that \,${{\bf U}_{n,k}^j=\ov{\bf P}_j\,{\bf U}_{n,k}}$,
\,where the projector \,${\ov{\bf P}_j:{\bf R}^{2d}\to{\bf R}^j}$ \,can be considered
as a linear operator with \,$\nnnorm{\ov{\bf P}_j}=1$ \,(see~(2.34)).
Applying Lemma~\ref {3.1} again, we obtain the relations
\,$\L({\bf U}_{n,k}^j)\in\AD{\sqrt 2\,\t}j$, \,$1\le j\le 2\4d$.

 It remains to verify that,
 for \,$ h\in{\bf R}^j$, \,${\norm h\sqrt 2\,\t<1}$,
\,the following inequality hold:
$$
\int\limits_{T}
\,\bbgl|\wh F_h(t)\bbgr|
\,dt\le\ffrac{(2\4\pi)^{j/2}\,\sqrt2\,\t\4j^{3/2}}
{\si\,(\det{\bf D})\ssqrt},\nopagebreak
\eqno(4.5)
$$
$$
T=\bgl\{t\in{\bf R}^j:4\4\norm t\4\sqrt 2\,\t\4j\ge1\bgr\},
\eqno(4.6)
$$
where
\,$F=\L \big({\bf U}_{n,k}^j\big)$,
\,and \,$\si^2$ \,is
the minimal eigenvalue of
~\,${{\bf D}=\cov {\bf U}_{n,k}^j}$. \,Note that,
according to~(3.9), we have
$$
{\bf D}=2^n\4{\bf I}_j,\quad
\si^2=2^n,\quad
\det{\bf D}=2^{nj}.
\eqno(4.7)
$$

Introduce $2^{n-1}$ random vectors
$$
{\bf X}_r=\big(X_r,\,X_{2^{n-1}+r}\big)\in{\bf R}^{2d},\qquad
r=2^{n-1}\cdot 2\4k+1,\dots,2^{n-1}\4(2\4 k+1).
\eqno(4.8)
$$
Obviously, these vectors are independent.
According to~(2.36),~(4.37) and (4.8),
$$
\wt{\bf U}_{n,k}
=\bgl(U_{n-1,2k},\,U_{n-1,2k+1})=
\sum_{r=2^{n-1}\cdot 2k+1}^{2^{n-1}(2k+1)}\3
{\bf X}_r.
\eqno(4.9)
$$

Denote now \,${\YY R{h}s}=\ovln{\L(X_s)}(h)$, \,for \,$s=1,\dots,
2^N$, \,${h\in\Rd}$, \,and \,${\YY M{h}r}\=\ovln{\L({\bf
X}_r)}(h)$, \,${{\YY Q{h}r}\=\ovln{\L({\bf A}\4{\bf X}_r)}(h)}$,
\,for
 \,$r=2^{n-1}\cdot 2\4k+1,\dots,2^{n-1}\4(2\4 k+1)$,
\,$h\in{\bf R}^{2d}$. \,As usually, we consider only such \,$h$
\,for which these distributions exist.
 Using~(2.8), we see that,
 for all \,$t\in{\bf R}^{2d}$,
\begin{eqnarray*}
{\wh Q_h^{(r)}(t)}
=\ffrac{\E \exp\bgl(\<h+i\4t,\4{\bf A} {\bf X}_r\>\bgr)}
{\E \exp\bgl(\<h,{\bf A} {\bf X}_r\>\bgr)}
\hskip-.6cm&&=\ffrac{\E \exp\bgl(\<{\bf A}^* h+i\4{\bf A}^* t,
 {\bf X}_r\>\bgr)}
{\E \exp\bgl(\<{\bf A}^* h, {\bf X}_r\>\bgr)}\\
&&=\,\wh M_{{\bf A}^*\! h}^{(r)}({\bf A}^* t).
\hbox{\rlap{\hskip3.14cm(4.10)}}
\end{eqnarray*}
By~(2.3) and~(4.9), we have (for \,$j=2\4d$)
$$
\bbgl|\wh F_h(t)\bbgr|
= \prod_{r=2^{n-1}\cdot 2k+1}^{2^{n-1}(2k+1)}\3\
\bbgl|\wh Q_h^{(r)}(t)\bbgr|.
\eqno(4.11)
$$
Split \,$t=\big(t_1,\dots,t_{2d}\big)\in{\bf R}^{2d}$ \,as
\,$t=\big(\YY t{}1,\,\YY t{}2\big)$, \,where we denote
 \,$\YY t{}1=\big(t_1,\dots,t_d\big)$ \,and
\,${\YY t{}2=\big(t_{d+1}},\dots,t_{2d}\big)\in\Rd$.
\,Using~formulae (2.8) and~(4.8) and introducing a similar
notation for~ \,$h\in{\bf R}^{2d}$, \,it is easy to check that
$$
\wh M_{h}^{(r)}(t)=\wh R_{\YY h{}1}^{(r)}\big(\YY t{}1\big)
\,\wh R_{\YY h{}2}^{(2^{n-1}+r)}\big(\YY t{}2\big).
\eqno(4.12)
$$
Note that
$$
\nnnorm {\2t\2}^2=
\nnnorm {\4t^{( 1)}}^2+
\nnnorm {\4t^{( 2)}}^2.
\eqno(4.13)
$$
\medbreak

\noindent{\it End of the proof of Theorem\/~{\rm \ref{1.1}}} \
Let now the distributions \,$\L(\xi_s)$
\,satisfy
the conditions of Theorem \ref {1.1}.
\,In this case, according to~(2.3), we have
\,${\YY R{h}s=
\ovln H_s(h)\,\ovln G(h)}$.
\,It is well-known that the conjugate distributions \,$\ovln G(h)$
\,of the Gaussian distribution~\,$G$
\,are also Gaussian with
covariance operator~\,$\cov \ovln G(h)=\cov G=b^2\,{\bf I}_d$.
\,Therefore,
$$
\bbgl|\wh R_h^{(s)}(t)\bbgr|
\le \exp\big(-b^2\nnnorm t^2\!/2\big),\qquad
t,h\in\Rd,\ \,\norm h\t<1.
\eqno(4.14)
$$
Using~(4.12)--(4.14),
we get, for \,$ t,h\in{\bf R}^{2d}$,  \,$\norm h\t<1$:
$$
\bbgl|\wh M_h^{(s)}(t)\bbgr|\le
\prod_{\mu=1}^2\exp\big(-b^2\nnnorm {\4t^{(\mu)}}^2\!/2\big)
=\exp\big(-b^2\nnnorm {\2t\2}^2\!/2\big).
\eqno(4.15)
$$
Applying
(2.40), (4.10) and (4.15)
with \,$t={\bf A}^*\4u$ \,and \,$h={\bf A}^*\4\gamma$,
\,we see that
$$
\bbgl|\wh Q_\gamma^{(s)}(u)\bbgr|\le
\exp\bgl(-b^2\4\nnorm {{\bf A}^*u}^2\!/2\bgr)
\le \exp\big(-b^2\nnnorm {u}^2\big),
\eqno(4.16)
$$
for \,$u, \gamma\in{\bf R}^{2d}$,  \,$\norm \gamma\sqrt 2\,\t<1$.
The relations~(4.11) and~(4.16) imply that
$$
\bbgl|\wh F_h(t)\bbgr|
\le \exp\bgl(-b^2\nnnorm t^2\cdot2^{n-2}\bgr),\qquad
t,h\in{\bf R}^j,\ \,\norm h\sqrt 2\,\t<1.
\eqno(4.17)
$$
It is clear that it suffices to verify~(4.17)
for
 \,$j=2\4d$ \,(for \,${1\le j<2\4d}$
\,one should apply~(4.17) for \,$j=2\4d$ \,and for
\,$t, h\in{\bf R}^{2d}$,
\,with \,$h_m=t_m=0$, \,${m=j+1,\dots,2\4d}$).

Using (4.6), (4.7) and (4.17),
 we see that
$$
\int\limits_{T}
\,\bbgl|\wh F_h(t)\bbgr|\,dt\le
\,\exp\Bigl(-\ffrac{b^2\cdot2^{n-3}}{32\,\t^2\4j^{2}}\Bigr)
\int\limits_{{\bf R}^j}
\,\exp\bgl(-b^2\4\nnnorm t^2\cdot2^{n-3}\bgr)\,dt\
$$
$$
{}\hskip-.4cm
= \ffrac {(2\4\pi)^{j/2}}{(b^2\cdot2^{n-2})^{j/2}}
\,\exp\Bigl(-\ffrac{b^2\cdot2^n}{2^8\,\t^2\4j^{2}}\Bigr)
$$
$$
{}\hskip.2cm\le\ffrac{(2\4\pi)^{j/2}\,\t^{4j\cdot2^n}}
{(\det{\bf D})\ssqrt\,\t^{2j}}
\le\ffrac{(2\4\pi)^{j/2}\,\t}
{2^{n/2}\, (\det{\bf D})\ssqrt},
\eqno(4.18)
$$
if \,$c_1$ \,is small enough. The relations (4.7) and (4.18)
imply~(4.5). It remains to apply Theorem~\ref {2.1} to complete
the proof of Theorem \ref {1.1}. $\square$\medbreak

\noindent{\it End of the proof of Theorem~{\rm\ref{1.2}}} \
Let now the distributions \,$\L(\xi_s)$
\,satisfy
the conditions of Theorem~\ref{1.2}.
\,In this case, according to~(2.8) and~(4.2), we have
\begin{eqnarray*}
\bbgl|\wh R_h^{(s)}(t)\bbgr| \hskip-.5cm&& = \,\biggl|
{\ffrac{\bgl(4+\t^2\4(d+\nnnorm h^2+2\4 i\8h, t\9-\nnnorm
t^2)\bgr) \,\exp\big((\nnnorm h^2+2\4 i\8h, t\9-\nnnorm
t^2)/2\big)} {\bgl(4+\t^2\4(d+\nnnorm h^2)\bgr) \,\exp\big(\nnnorm
h^2\!/2\big)}}\biggr|
\\
&&\le \,\bgl(2+\nnnorm t^2\big)\,\exp\big(-\nnnorm t^2\!/2\big)
\\
&&\le \,c_{25}\,\exp\big(-\nnnorm t^2\!/4\big),\qquad\qquad \norm
h\t<1. \hbox{\rlap{\hskip5.48cm(4.19)}}
\end{eqnarray*}
The rest of the proof is omitted. It is similar to that of
Theorem \ref {1.1}  with~\,${b^2=\half}$. \,The presence of~\,$c_{25}$
\,in the right-hand side of~(4.19) can be easily compensated
by choosing \,$c_1$ \,to be sufficiently small.
\medbreak

\noindent{\it End of the proof of Theorem\/~{\rm \ref{1.3}}} \
Consider the dyadic scheme with
$$
\L(\xi_s)=\L(X_s)=\YY L{}s\4P, \qquad s=1,\dots, 2^N.
\eqno(4.20)
$$

Putting \,$H\=\L(\zeta)$,
\,\,${\psi_h(x)=e^{\8h, x\9}\,p(x)}$,
\,$h,x\in{\bf R}^d$, \,and
integrating by parts, we see that
(for~\,$t\in{\bf R}^d$, \,$t\ne0$)
\begin{eqnarray*}
\wh H_h(t)&=&\bgl(\!\E e^{\8h,\zeta\9}\bgr)^{-1}
\int\limits_{\nnnorm x\le b_1}
e^{i\2\8t,x\9}\,\psi_h(x)\,dx\\
&=&-\bgl(\!\E e^{\8h,\zeta\9}\bgr)^{-1} \int\limits_{\nnnorm x\le
b_1} \ffrac{e^{i\2\8t,x\9}} {i\4\nnnorm t^2}d_t \4\psi_h(x)\,dx,
\hbox{\rlap{\hskip2.6cm(4.21)}}
\end{eqnarray*}
where \,$H_h=\ovln H(h)$.
\,Besides, using~(1.9), we see that
$$
\sup_{\nnnorm x\le b_1}\,\sup_{\nnnorm h \4b_2\le1}
\,\bgl|d_t \4\psi_h(x)\bgr|\le b_5\4\nnnorm t.
\eqno(4.22)
$$
As in the formulation of Theorem~\ref {1.3}  we denote
by \,$b_m$ \,different positive
quantities depending on~\,$H$. \,Note
that the quantities depending on the dimension~\,$d$
\,can be considered as depending on~\,$H$ \,only as well.
From (4.21) and
~(4.22) it follows that
$$
\sup_{\nnnorm h \4b_2\le1}\,\bbgl|\wh H_h(t)\bbgr|
\le b_6\4\nnnorm t^{-1}
\eqno(4.23)
$$
(note that, by the Jensen inequality,
$\E e^{\8h,\zeta\9}\ge e^{\E\8h,\zeta\9}=1$).
The inequality~(4.23) implies that
$$
\sup_{\nnnorm h \4b_2\le1}
\,\bbgl|\wh H_h(t)\bbgr|\le
 \Big(\,1+\ffrac{\nnnorm t}{b_7}\Big)^{\!-1}
\qquad\hbox{for}\quad \nnnorm t\ge b_7=2\4b_6
\nopagebreak
\eqno(4.24)
$$
and
$$
\sup_{\nnnorm h \4b_2\le1}\,\sup_{\nnnorm t\ge b_7}
\,\bbgl|\wh H_h(t)\bbgr|\le \half.\nopagebreak
\eqno(4.25)
$$

Since the distributions \,$H_h$ \,are absolutely continuous,
the relation \,$\bbgl|\wh H_h(t)\bbgr|=1$
\,can be valid for \,$t=0$ \,only. Furthermore,
the function \,$\bbgl|\wh H_h(t)\bbgr|$ \,considered
as a function of two variables \,$h$ \,and \,$t$ \,is
continuous for all \,$h,t\in{\bf R}^d$. \,Therefore,
$$
\sup_{\nnnorm h \4b_2\le1}
\,\sup_{b_8\le\4\nnnorm t\4\le\4 b_7}
\,\bbgl|\wh H_h(t)\bbgr|\le b_9<1,\nopagebreak
\eqno(4.26)
$$
where
$$
b_8=\big(4\sqrt2\,b_2\4d\big)\me\quad\hbox{and}
\quad b_9\ge\half.
\eqno(4.27)
$$
The inequalities (4.25) and (4.26) imply that
$$
\sup_{\nnnorm h \4b_2\le1}
\,\sup_{\nnnorm t\4\ge\4 b_8}
\,\bbgl|\wh H_h(t)\bbgr|\le b_9\=e^{-b_{10}}<1.
\eqno(4.28)
$$
Denoting \,$\YY Lhs=\ovln L{}^{(s)}(h)$,
\,$h\in \Rd$, \,$s=1,\dots,2^N$,
\,and using~(1.11),~(1.12), (2.3) and~(2.8),
it is easy to see that
$$
\wh R{}^{(s)}_h(t)=\bgl(\wh H_{h/\sqrt m}\big(t/\sqrt m\big)\bgr)^{\!m}
\,\wh L{}^{(s)}_h(t).
\eqno(4.29)
$$
The relations~(1.10), (4.24), (4.28) and~(4.29) imply that
$$
\sup_{\nnnorm h \4\t\le1}
\,\bbgl|\wh R{}^{(s)}_h(t)\bbgr|\le
 \Big(\,1+\ffrac{\nnnorm t}{b_7\,\sqrt m}\Big)^{\!-m}\for
\nnnorm t\ge b_7\,\sqrt m
\eqno(4.30)
$$
and
$$
\sup_{\nnnorm h \4\t\le1}
\,\sup_{\nnnorm t\4\ge\4 b_8\sqrt m}
\,\bbgl|\wh R{}^{(s)}_h(t)\bbgr|\le e^{-mb_{10}}.
\eqno(4.31)
$$

Using~(4.12), (4.13), (4.20) and (4.30),
we get, for
 \,$r=2^{n-1}\cdot 2\4k+1,\dots,2^{n-1}\4(2\4 k+1)$,
\,$\nnnorm t\ge b_7\,\sqrt {2\4m}$,
\,$t\in{\bf R}^{2d}$,
$$
\sup_{\htau\le 1}
\,\bbgl|\wh M_h^{(r)}(t)\bbgr|\le
\min_{\mu=1,2}
\Big(\,1+\ffrac{\nnnorm {t^{(\mu)}}}{b_7\,\sqrt m}\Big)^{\!-m}
\le \Big(\,1+\ffrac{\nnnorm t}{b_7\,\sqrt{2\4m}}\Big)^{\!-m}.
\eqno(4.32)
$$
Moreover,
$$
\sup_{\nnnorm h \4\t\le1}
\,\sup_{\nnnorm t\4\ge\4 b_8\sqrt {2\4m}}
\,\bbgl|\wh M_h^{(r)}(t)\bbgr|\le e^{-mb_{10}}.
\eqno(4.33)
$$
 Using~(2.40),~(4.10), (4.32) and~(4.33), we see that,
for the same~\,$r$
\,and for \,$t\in{\bf R}^{2d}$,
\,$\nnnorm t\ge b_7\,\sqrt m$,
$$
\sup_{\htau\sqrt2\le 1}
\,\bbgl|\wh Q_h^{(r)}(t)\bbgr|\le
\Big(\,1+\ffrac{\nnnorm t}{b_7\,\sqrt m}\Big)^{\!-m}
\eqno(4.34)
$$
and
$$
\sup_{\htau\sqrt2\le 1}
\,\sup_{\nnnorm t\4\ge\4 b_8\sqrt m}
\,\bbgl|\wh Q_h^{(r)}(t)\bbgr|\le e^{-mb_{10}}.
\eqno(4.35)
$$

It is easy to see that the relations~(4.11),~(4.34) and~(4.35) imply that,
 for \,$ h\in{\bf R}^j$, \,${\norm h\sqrt 2\,\t<1}$,
\,and for \,$t\in{\bf R}^j$,
\,$\nnnorm t\ge b_7\,\sqrt m$,
$$
\bbgl|\wh F_h(t)\bbgr|\le
\Big(\,1+\ffrac{\nnnorm t}{b_7\,\sqrt m}\Big)^{\!-m\cdot2^{n-1}}
\eqno(4.36)
$$
and
$$
\sup_{\nnnorm t\4\ge\4 b_8\sqrt m}
\,\bbgl|\wh F_h(t)\bbgr|\le e^{-mb_{10}\cdot2^{n-1}}.
\eqno(4.37)
$$
It suffices to prove (4.36) and~(4.37)
for \,$j=2\4d$ \,(for \,$1\le j<2\4d$
\,one should apply (4.36) and~(4.37) for \,$j=2\4d$ \,and for
\,$ h\in{\bf R}^{2d}$, \,${\norm h\sqrt 2\,\t<1}$, \,$t\in{\bf R}^{2d}$
\,with \,$h_m=t_m=0$, \,$m=j+1,\dots,2\4d$).

Note now that the set~\,$T$ \,defined in~(4.6)
satisfies the relation
$$
T\subset\bgl\{t\in{\bf R}^j:\nnnorm t\4\ge\4 b_8\4\sqrt m\bgr\}
\eqno(4.38)
$$
(see~(1.10) and~(4.27)). Below (in the proof of (4.5)) we assume
that ${\norm h\sqrt 2\,\t<1}$. \,According to~(4.37) and~(4.38),
for \,$t\in T$ \,we have
$$
\bbgl|\wh F_h(t)\bbgr|\ssqrt\le e^{-m b_{10}\cdot2^{n-2}}.
\eqno(4.39)
$$
Taking into account that \,$\bbgl|\wh F_h(t)\bbgr|\le1$,
 \,and \,$m\ge b_4$, \,choosing~\,$b_4$ \,to be
sufficiently large
and using (1.10), (4.7), (4.36) and~(4.39), we obtain
\begin{eqnarray*}
\int\limits_{T} \,\bbgl|\wh F_h(t)\bbgr|\,dt \,\5&\le&\5
\,\exp\big(-m\4 b_{10}\cdot2^{n-2}\big) \biggl(\,\int\limits_{{\bf
R}^j} \,\Big(\,1+\ffrac{\nnnorm t}{b_7\,\sqrt
m}\Big)^{\!-m\cdot2^{n-2}} dt+b_{11}\4m^{d/2}\biggr)\\
\vspace{.5pt} \5&\le&\5 \,b_{12}\,m^{d/2} \,\exp\big(-m\4
b_{10}\cdot2^{n-2}\big)\\ \vspace{1\jot} \5&\le&\5\ffrac
{(2\4\pi)^{j/2}\sqrt2\,b_2\4j^{3/2}} {m\ssqrt\cdot2^{n/2}\cdot2^{n
j/2}}= \ffrac{(2\4\pi)^{j/2}\,\sqrt2\,\t\4j^{3/2}} {\si\,(\det{\bf
D})\ssqrt}. \hbox{\rlap{\hskip4.03cm(4.40)}}
\end{eqnarray*}
The inequality (4.5) follows from~(4.40) immediately. It remains
to apply Theorem~\ref {2.1}. ~$\square$\medbreak

\noindent{\it Proof of Theorem\/~{\rm \ref{1.4}}} \
Define \,$m_0,m_1,m_2,\dots$ \,and
\,$n_1,n_2,\dots$ \,by
$$
m_0=0,\quad m_s=2^{2^s}, \qquad
n_s=m_s-m_{s-1},\qquad s=1,2,\dots.
\eqno(4.41)
$$
It is easy to see that
$$
\log n_s\le\log m_s=2^s\4\log2,\qquad s=1,2,\dots.
\eqno(4.42)
$$
By Assertion~A (see~(1.5)), for any \,$s=1,2,\dots$
\,one can construct on a probability space
a sequence of i.i.d.\  \,$\YY X1s,\dots,\YY X {n_s}s$
\,and a sequence of i.i.d.\  Gaussian
\,$\YY Y1s,\dots,\YY Y {n_s}s$ \,so that
\,$\L(\YY Xks)=\L(\xi)$, \,$ \E \YY Yks=0$,
\,$\cov \YY Yks={\bf I}_d$,
\,and
$$
{}\P\bgl\{\4c_2\,\DE_s\ge \t\4d^{3/2}\bgl(
c_3\,\log^*d\,\log n_s
+x\bgr)\4\bgr\}
\le e^{-x},\qquad x\ge0,
\eqno(4.43)
$$
where
$$
\DE_s=\max_{1\le r\le n_s}
\,\Bigl|\,\sum\limits_{k=1}^r X_k^{(s)}
-\sum\limits_{k=1}^r Y_k^{(s)}\,\Bigr|.
\eqno(4.44)
$$

It is clear that we can define all the vectors mentioned above
on the same probability space so that the collections
\,$\Xi_s=\bgl\{\YY X1s,\dots,\YY X {n_s}s;
\,\YY Y1s,\dots,\YY Y {n_s}s\bgr\}$, \,$s=1,2,\dots$
\,are jointly independent. Then we define
\,$X_1,X_2,\dots$ \,and
\,$Y_1,Y_2,\dots$ \,by
$$
\begin{tabular}{r c l}
$X_{m_{s-1}+k}\5$ &=& $\5X_k^{(s)},$\\
$Y_{m_{s-1}+k}\5$ &=& $\5Y_k^{(s)},$
\end{tabular}
\qquad k=1,\dots,n_s,\quad s=1,2,\dots.
\eqno(4.45)
$$
In order to show that these sequences satisfy the assertion
of Theorem \ref {1.4}, it remains to verify
the equality~(1.13).

Put
$$
c_{25}=
\ffrac{(c_3\4\log2+1)}
{c_2},\qquad
c_{26}=c_{25}\sum_{l=0}^\infty2^{-l/2}=\ffrac{c_{25}\4
\sqrt2}
{\sqrt2-1},
\eqno(4.46)
$$
and introduce the events
$$
A_l=\bgl\{\4\om:\DE^{(l)}\ge
2^l\,c_{26}\,\t\4d^{3/2}\log^*d\4\bgr\},
\qquad l=1,2,\dots,
\eqno(4.47)
$$
where
$$
\DE^{(l)}=\max_{1\le r\le m_l}
\,\Bigl|\,\sum\limits_{j=1}^r X_j
-\sum\limits_{j=1}^r Y_j\,\Bigr|.
\eqno(4.48)
$$
According to (4.44), (4.45) and~(4.48), we have
$$
\DE^{(l)}\le \DE_1+\dots+\DE_l.
\eqno(4.49)
$$
Taking into account the relations
~(4.42), (4.46),~(4.47),~(4.49) and applying the inequality~(4.43)
 with~\,$x=2^{(s+l)/2}$,
\,we get
\begin{eqnarray*}
{}\P\bgl\{A_l\bgr\}&\le &\sum_{s=1}^l
\P\bgl\{\4\DE_s\ge2^{(s+l)/2}\,c_{25}\,\t\4d^{3/2}\log^*d\4\bgr\}\\
&\le& \sum_{s=1}^l\exp\bgl(-2^{(s+l)/2}\bgr) \le
c\,\exp\bgl(-2^{l/2}\bgr). \hbox{\rlap{\hskip2.3cm(4.50)}}
\end{eqnarray*}
The inequality~(4.50) implies that \,$\sum\limits_{l=1}^\infty
\P\bgl\{A_l\bgr\}<\infty $, \,Hence, by the Borel--Cantelli lemma
with probability one a finite number of the events~\,$A_l$
\,occurs only. This implies the equality~(1.13) with
\,$c_4=2\4c_{26}\big/\log2$ \,(see~(4.41), (4.47) and~(4.48)).
$\square$\medbreak

\vspace*{30pt}

\section*{References}
\begin{enumerate}

\item
 B\'artfai, P. (1966).
 Die Bestimmung der zu einem wiederkehrenden
Prozess geh\"orenden Verteilungfunktion aus den mit
Fehlern behafteten Daten einer einzigen Realisation,
{\it Studia Sci. Math. Hungar.}, {\bf  1},
 161--168.

\item
 Dudley, R. M. (1989).
{\it Real analysis and probability},
 Pacific Grove, California: Wads\-worth \& Brooks/Cole.

\item
 Einmahl, U. (1989).
 Extensions of results of Koml\'os, Major and Tusn\'ady
to the multivariate case,
{\it J.~Multivar. Anal.}, {\bf  28},  20--68.

\item
 Etemadi, N. (1985).
 On some classical results in probability theory,
{\it San\-khy$\bar{\hbox{a}}$, Ser.~A}, {\bf  47}, 2,
  215--221.
\item
 Hoffmann-J\o rgensen, J. (1994).
{\it Probability with a view toward statistics},~\rm I,
 New York: Chapman \& Hall.

\item
 G\"otze, F.
and  Zaitsev A. Yu. (1997).
 Multidimensional
Hungarian construction for vectors
with almost Gaussian smooth distributions,
 {\it Preprint 97-071 SFB~343}, Universit\"at Bielefeld.

\item
 Koml\'os, J., Major, P., Tusn\'ady, G. (1975-76).
 An approximation of partial sums of independent RV'-s
and the sample DF.\/ {\rm  I; II},
{\it Z. Wahrscheinlichkeitstheor. verw. Geb.}, {\bf  32},   111--131;
 {\bf  34},  34--58.

\item
 Major, P. (1978).
 On the invariance principle for sums of independent
identically distributed random variables
{\it J. Multivar. Anal.}, {\bf  8},
 487--517.

\item
 Massart, P. (1989).
 Strong approximation for multivariate
empirical and related
processes, via KMT construction,
{\it Ann. Probab.}, {\bf  17}, 1,
 266--291.

\item
 Rosenblatt, M. \!(1952).
 Remarks on a multivariate transformation,
{\it Ann. Math. Statist.}, {\bf  23},
 470--472.

\item
 Sakhanenko, A. I. (1984).
 Rate of convergence in the invariance principles
for variables with exponential moments that are not
identically distributed,
 In: {\it Trudy Inst. Mat. SO AN SSSR}, {\bf  3}, pp.  4--49,
 Novosibirsk: Nauka (in Russian).

\item
 Zaitsev, A. Yu. (1986).
 Estimates of the L\'evy--Prokhorov distance in the
 multivariate central limit theorem for random variables with finite
 exponential moments,
 {\it Theor. Probab. Appl.}, {\bf  31},  2,
  203--220.

\item
 Zaitsev, A. Yu. (1995).
 Multidimensional version of the results
 of Koml\'os, Major and Tusn\'ady for
 vectors with finite exponential moments,
 {\it Preprint 95-055 SFB~343}, Universit\"at Bielefeld.

\item
 Zaitsev, A. Yu. (1996).
 Estimates for quantiles of smooth conditional
distributions and multidimensional invariance principle,
 {\it Siberian Math.~J.}, {\bf  37},  4,
  807--831
(in Russian).

\item
 Zaitsev, A. Yu. (1998a).
 Multidimensional version of the results
of Koml\'os, Major and Tusn\'ady for
vectors with finite exponential moments,
 {\it  ESAIM : Probability and Statistics}, {\bf 2}, 41--108.

\item
 Zaitsev, A. Yu. (1998b).
 Multidimensional version of the results
 of Sakhanenko in the invariance principle for
 vectors with finite exponential moments,
 {\it Preprint 98-045 SFB~343}, Universit\"at Bielefeld.

\end{enumerate}

\end{document}